\documentclass[review,12pt]{elsarticle}

\usepackage{setspace}
\usepackage{subcaption}
\usepackage{multirow}
\usepackage{multicol}
\usepackage{graphicx}
\RequirePackage{amsmath}
\usepackage{bm}
\usepackage{xfrac}
\usepackage{float}
\usepackage{arydshln} 

\usepackage{amsmath,esint,amsfonts}

\usepackage[most]{tcolorbox}
\usepackage[colorinlistoftodos]{todonotes}

\usepackage{lineno}
\usepackage[colorlinks=true,linkcolor=blue,urlcolor=blue,citecolor=blue,anchorcolor=blue]{hyperref}
\usepackage{cleveref} 

\journal{Journal}

\biboptions{sort&compress} 

\begin{document}
\begin{frontmatter}

	\title{Consistency and Convergence of a High Order Accurate Meshless Method for Solution of Incompressible Fluid Flows}

	\author{Shantanu Shahane\fnref{Corresponding Author}$^{a,}$}
	\author{Surya Pratap Vanka$^{b}$}
	\address{a National Center for Supercomputing Applications,\\
		University of Illinois at Urbana-Champaign, Urbana, Illinois 61801}
	\address{b Department of Mechanical Science and Engineering,\\
		University of Illinois at Urbana-Champaign, Urbana, Illinois 61801 \vspace{-0.7cm}}
	\fntext[Corresponding Author]{\vspace{0.3cm}Corresponding Author Email Address: \url{shahaneshantanu@gmail.com}}

	\begin{abstract}
		Computations of incompressible flows with all velocity boundary conditions require solution of a Poisson equation for pressure or pressure-corrections with all Neumann boundary conditions. Discretization of such a Poisson equation results in a rank-deficient ill-conditioned matrix of coefficients. When a non-conservative discretization method such as finite difference, finite element, or spectral scheme is used, the ill-conditioned matrix also generates an inconsistency which makes the residuals in the iterative solution to saturate at a threshold level that depends on the spatial resolution and the order of the discretization scheme. In this paper, we examine this inconsistency for a high-order meshless discretization scheme suitable for solving the equations on a complex domain. The high order meshless method uses polyharmonic spline radial basis functions (PHS-RBF) with appended polynomials to interpolate scattered data and constructs the discrete equations by collocation. The PHS-RBF provides the flexibility to vary the order of discretization by increasing the degree of the appended polynomial. In this study, we examine the convergence of the inconsistency for different spatial resolutions and for different degrees of the appended polynomials by solving the Poisson equation for a manufactured solution as well as the Navier-Stokes equations for several fluid flows. We observe that the inconsistency decreases faster than the error in the final solution, and eventually becomes vanishing small at sufficient spatial resolution. The rate of convergence of the inconsistency is observed to be similar or better than the rate of convergence of the discretization errors. This beneficial observation makes it unnecessary to regularize the Poisson equation by fixing either the mean pressure or pressure at an arbitrary point. A simple point solver such as the SOR is seen to be well-convergent, although it can be further accelerated using multilevel methods.\\
	\end{abstract}


\end{frontmatter}

\section{Introduction}

Accurate computations of incompressible flows in complex domains are required to understand transport processes in a variety of engineering flow devices. For an incompressible flow, the Mach number (ratio of a characteristic velocity to sound speed) is low, and the fluid density is a weak function of the local pressure. Hence, it is commonly assumed that the density at a point is either constant or only a function of variables such as temperature and composition. The local pressure variations only drive the velocity field without affecting the local density. Since there is no explicit equation for pressure, the pressure field is determined such that the velocity field is divergence-free. Computational techniques to solve incompressible flows are therefore labeled as pressure based methods versus the term density based often used for computing compressible flows \cite{anderson2016computational, ferziger2002computational, roache1998fundamentals}. Two main types of pressure based methods have been developed in the literature. The first type, called pressure projection methods \cite{harlow1965numerical, chorin1967numerical, kim1985application}, solves the momentum equations for an intermediate velocity field without the pressure gradient term. This intermediate velocity is then projected to a divergence-free state by solving a pressure Poisson equation that has the divergence in the intermediate velocity field as its source term. Pressure projection methods march in time by solving the momentum equations followed by the pressure Poisson equation and corrections to the intermediate velocity fields. The pressure projection concept has been used in conjunction with finite-difference, finite element, finite volume, and spectral methods \cite{mittal2021multirate, rosales2021high, shahane2019finite}. Further, methods with both explicit and implicit formulations of the diffusion terms have been developed. A second approach is a semi-implicit formulation in which the pressure field is determined through an iterative correction strategy. This approach is followed in algorithms originated from SIMPLE \cite{patankar1983calculation, patankar2018numerical, issa1986computation, shahane2021numerical} in which a pressure correction equation is derived from truncated momentum equations and the discrete continuity equation. The pressure corrections eventually become small when the correct pressure and velocity fields converge to satisfy momentum and continuity equations.
\par While both pressure projection and pressure correction techniques approach the same solution and can be shown to be somewhat equivalent, the convergence rates of the two methods can be different because of their respective formulations. Further, the semi-implicit formulation was developed primarily with the steady-state equations in mind while the pressure projection algorithms were focused on time-evolving flows. However, despite this difference in the two techniques, both assume that density is not a function of local pressure but only the system pressure. This has an important consequence on computation of incompressible flows with purely velocity boundary conditions. With velocity boundary conditions, the pressure field can only be determined to an arbitrary constant \cite{patankar2018numerical} because only the pressure differences are important and not the actual pressures. Any uniform arbitrary constant added to the calculated pressure is also a solution to the problem and gives the same velocity field.
\par Another issue with the solution of the pressure Poisson (or pressure correction) equation is the consistency of the discrete equations. One requirement of incompressible flow computations is satisfaction of the global continuity equation (sum of inflows equals sum of outflows). This can be clearly ensured in finite volume methods where mass fluxes are defined on the faces of the finite volumes. When the discrete mass imbalances (or divergences in the intermediate velocities) are calculated, mass fluxes are first evaluated on the cell faces and summed to calculate local divergences. Thus, in finite volume methods, the sum of local divergences is equal to the global divergences, which becomes zero since the boundary velocities satisfy global balance. Similarly, in the discretization of the pressure Poisson equation
\begin{equation}
\nabla \bullet (\nabla p) = \frac{\rho}{\Delta t} \left(\frac{\partial \hat{u}}{\partial x} + \frac{\partial \hat{v}}{\partial y} \right)
\end{equation}
where, the $\nabla p$ is evaluated at the faces of the finite volumes, with the result that the sum of the coefficients of the discrete pressures also becomes zero. The pressure Poisson equation therefore becomes ill-conditioned and cannot be solved using direct solvers.
\par Two strategies have been proposed to overcome this difficulty. One approach is to replace one of the discrete equations with an equation which fixes the level of the pressure at that location. Thus, a new equation
\begin{equation}
p_{\text{iref}}=c
\end{equation}
is used to replace the discrete equation at cell `iref' which is an arbitrary location at which pressure is fixed to a value `$c$'. This strategy provides a consistent and convergent solution to the pressure Poisson equation with direct or iterative solutions. A second approach, commonly followed in both pressure correction as well as pressure projection algorithms is to use an iterative solver such as SOR or alternate-direction line inversions until convergence, and then subtract the calculated pressure at one reference location from the other discrete pressures. This fixes the pressure at the reference location to be zero after obtaining a converged solution. Both methods achieve the same result, but the convergence of the second method has been observed to be better.
\par When the pressure Poisson equation is discretized by any method other than the finite volume method, conservation principles are not used. For example, in finite difference methods, the derivatives at the grid nodes are calculated using a stencil derived from Taylor series expansions. There is no enforcement of mass balances. As a result, an inconsistency arises because the coefficients of the discrete pressures as well as the mass sources on the right-hand sides sum to different non-zero values. The numerical condition number of the coefficient matrix is also very high, since it is rank-deficient because of the arbitrariness of the pressure level. If an iterative method is used for the solution, the equations cannot be expected to converge below the inconsistency error, resulting in a `stationary' residual, whose value depends on the grid resolution and the order of discretization accuracy. For these non-conservative discretization methods also, fixing the pressure level at a reference location is possible either prior to solution or after convergence
\par Finite volume methods are limited by being second order accuracy unless complex reconstruction schemes are used. When high order reconstruction schemes are used for complex geometries discretized with unstructured polygonal elements, much care is needed to ensure that reciprocity is obeyed at the cell faces so that global mass balances and the cancellation of the discrete pressure coefficients happens. Otherwise, these too will become inconsistent akin to finite difference and finite element methods. In finite difference methods, the stationary level to which the pressure equation converges is comparable to the discretization accuracy and decreases with grid refinement as per the order of the discretization scheme. This phenomenon is encouraging because it is then not necessary to supplement the pressure Poisson equation with additional equations to fix either the mean or a local pressure.
\par There are several types of meshless methods used for the solutions of heat transfer and fluid flow problems \cite{shu2003local, ding2006numerical, sanyasiraju2008local, vidal2016direct, zamolo2019solution, bayona2017role_II, kosec2008solution, kosec2020radial, unnikrishnan2022shear, shankar2018hyperviscosity, bartwal2021application, shahane2021high, shahane2021semi}. The intent of this work is to study the consistency and convergence characteristics of a recently developed high order meshless technique \cite{shahane2021high} that solves the incompressible Navier-Stokes equations using Polyharmonic Spline Radial Basis Function (PHS-RBF) interpolations of scattered data. When the PHS-RBF are appended with polynomials, the discretization errors decrease as per the degree of the polynomial (with order $\mathcal{O}(k-1)$ for the Laplacian estimation, where $k$ is the degree of the polynomial). However, as in the case of the finite difference or finite element methods, the pressure Poisson equation also has consistency issues, and does not converge to round-off errors. The pressure iterations reach a stationary level which depends on an average spacing between the scattered points and the degree of the appended polynomial. In this paper, we study the iterative and discretization convergences of this high order meshless method in several flows in complex domains. We observe that these consistency errors have similar or higher order of convergence as the solution errors. All the research presented here is performed using the open source Meshless Multi-Physics Software (MeMPhyS) \cite{shahanememphys}.
\par \Cref{Sec:Solution of Poisson Equation with Manufactured Solution using Structured Grid Finite Difference and Volume Methods} illustrates the inconsistency issue in more detail by solving the two-dimensional Poisson equation with finite difference and finite volume formulations on a structured uniform Cartesian point placement and finite volumes respectively. We present the convergence rates for both the residuals and the L1 norm of the solution errors (difference between the exact solution and the numerical solution). It is seen that the orders of convergence with mesh size of both these errors are similar, leading to the conclusion that termination of the iterations when the residuals reach an asymptotic value does not introduce errors much larger than the discretization errors themselves. \Cref{Sec:Numerical Method} describes the high order meshless discretization on scattered points and a fractional-step algorithm utilizing the concept for solving incompressible Navier-Stokes equations. In \cref{Sec:Poisson Equation with Manufactured Solution Meshless Discretization}, we first present results of convergence for the Poisson equation with Neumann boundary conditions and a manufactured solution using the meshless method. Here we have solved the discrete equations using the inexpensive SOR algorithm. \Cref{Sec:Incompressible Fluid Flow Problems Manufactured Solution} considers two model flow problems for which exact solutions are available: the Kovasznay flow \cite{kovasznay1948laminar} and a circular Couette flow. In \cref{Sec:Incompressible Fluid Flow Problems Analysis of Complex Flows}, we consider two more fluid flow problems, in complex geometries, and systematically study consistency and convergence for different point spacings (h-refinement) and polynomials (p-refinement). The paper is concluded with main observations and plans for future work.

\section{Solution of Poisson Equation with Manufactured Solution using Structured Grid Finite Difference and Volume Methods} \label{Sec:Solution of Poisson Equation with Manufactured Solution using Structured Grid Finite Difference and Volume Methods}
\subsection{Origin of Inconsistency}
The inconsistency mentioned in this paper arises in all formulations that solve the incompressible Navier-Stokes equations with prescribed velocity conditions on all the boundaries. It does not arise if pressure is prescribed at any boundary, even on a partial segment. For a prescribed velocity boundary, the boundary pressure cannot be fixed. Rather the normal gradient is to be prescribed either to be zero or from the normal momentum equation. If all the boundaries are prescribed with velocity conditions, the boundary conditions for pressure become Neumann type. In such a case, there is no unique solution for the pressure equation as the pressure field can only be determined to an arbitrary constant (i.e for a given $p(\bm{x},t)$ field, $p(\bm{x},t) +c$ is also a solution). In the finite volume method, this property manifests in the form of one less discrete equation, as the sum of all the discrete pressure Poisson equations strictly sum to a null equation. Consider a one-dimensional Poisson equation discretized by a finite volume method. We consider a simple manufactured solution given by:
\begin{equation}
	p_{exact}(x,y)=\sin(\omega \pi x)
\end{equation}
where, $\omega$ is the number of waves in a domain of unit length for which, the Poisson equation is given as
\begin{equation}
	\frac{\partial^2 p}{\partial x^2} = -\omega^2 \pi^2 \sin(\omega \pi x)
	\label{Eq:PPE 1D}
\end{equation}
In the finite volume method, the discrete form of the partial differential equation is obtained by integrating the above equation over control volumes as
\begin{equation}
	\frac{\partial p}{\partial x}\biggr\rvert_{i+\frac{1}{2}} - \frac{\partial p}{\partial x}\biggr\rvert_{i-\frac{1}{2}} = \omega \pi \left[ \cos(\omega \pi x)\big|_{i+\frac{1}{2}} - \cos(\omega \pi x)\big|_{i-\frac{1}{2}} \right]
\end{equation}
where, the subscript indices $i$, $i+\frac{1}{2}$ and $i-\frac{1}{2}$ denote the center of a one dimensional control volume and its east and west face center locations respectively as shown in \cref{Fig:FVM 1D schematic}. The derivatives at face centers can be evaluated based on the cell centered values using the second ordered central difference method:
\begin{equation}
	\frac{p_{i+1}-p_i}{\Delta x} - \frac{p_i-p_{i-1}}{\Delta x} = \omega \pi \left[ \cos(\omega \pi x)\big|_{i+\frac{1}{2}} - \cos(\omega \pi x)\big|_{i-\frac{1}{2}} \right]
	\label{Eq:PPE 1D FVM}
\end{equation}
For demonstration, consider a one dimensional domain discretized by 6 control volumes:
\begin{figure}[H]
	\centering
	\includegraphics[width=0.6\textwidth]{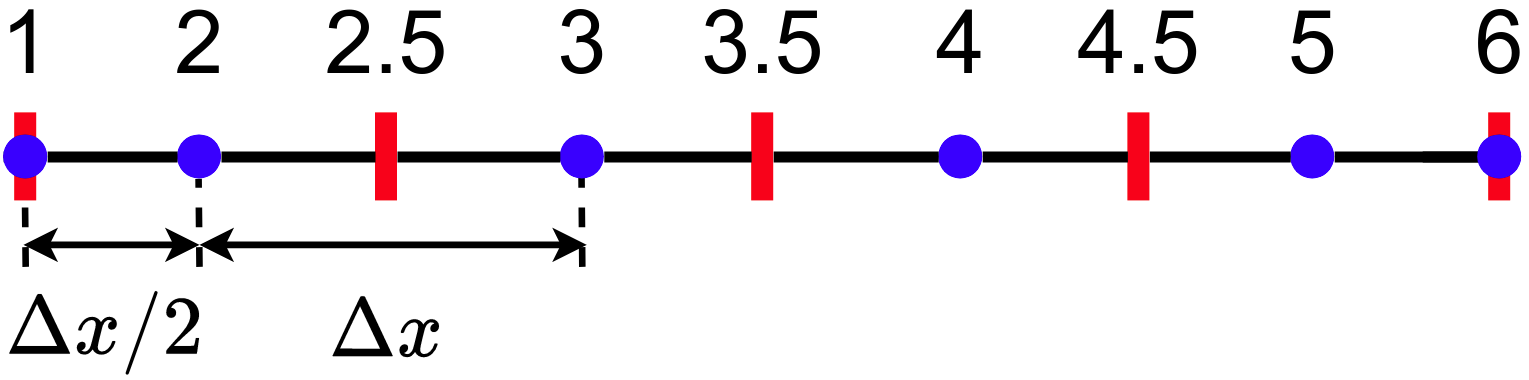}
	\caption{Schematic of Finite Volume Discretization}
	\label{Fig:FVM 1D schematic}
\end{figure}
On the boundaries, face locations 1.5 and 5.5 coincide with cell centers 1 and 6 respectively. Assembling \cref{Eq:PPE 1D FVM} with a Neumann boundary conditions at cells 1 and 6 gives:
\begin{equation}
	\begin{bmatrix}
		-2 & 2 & 0 & 0 & 0 & 0   \\
		2 & -3 & 1 & 0 & 0 & 0   \\
		0 & 1 & -2 & 1 & 0 & 0    \\
		0 & 0 & 1 & -2 & 1 & 0    \\
		0 & 0 & 0 & 1 & -3 & 2   \\
		0 & 0 & 0 & 0 & 2 & -2   \\
	\end{bmatrix}
	\begin{bmatrix}
		p_1  \\
		p_2  \\
		p_3  \\
		p_4  \\
		p_5  \\
		p_6  \\
	\end{bmatrix} =
	\begin{bmatrix}
		[S(x_1)] \Delta x  \\
		[S(x_{2.5}) -S(x_1)] \Delta x  \\
		[S(x_{3.5}) -S(x_{2.5})] \Delta x  \\
		[S(x_{4.5}) -S(x_{3.5})] \Delta x  \\
		[S(x_{6}) -S(x_{4.5})] \Delta x  \\
		[-S(x_{6})] \Delta x  \\
	\end{bmatrix}
	\label{Eq:PPE 1D FVM matrix}
\end{equation}
where, $S(x)=\omega \pi \cos(\omega \pi x)$, the source term of \cref{Eq:PPE 1D FVM}. We can see that when the individual equations are added, both the left and right sides sum to zero. On the other hand, when the \cref{Eq:PPE 1D} is discretized over a uniform grid of 6 points using the finite difference method (FDM), we get:
\begin{equation}
	\begin{bmatrix}
		-1 & 1 & 0 & 0 & 0 & 0   \\
		1 & -2 & 1 & 0 & 0 & 0   \\
		0 & 1 & -2 & 1 & 0 & 0    \\
		0 & 0 & 1 & -2 & 1 & 0    \\
		0 & 0 & 0 & 1 & -2 & 1   \\
		0 & 0 & 0 & 0 & 1 & -1   \\
	\end{bmatrix}
	\begin{bmatrix}
		p_1  \\
		p_2  \\
		p_3  \\
		p_4  \\
		p_5  \\
		p_6  \\
	\end{bmatrix} =
	\begin{bmatrix}
		\omega \pi \cos(\omega \pi x_1) \Delta x  \\
		-\omega^2 \pi^2 \sin(\omega \pi x_2) \Delta x^2  \\
		-\omega^2 \pi^2 \sin(\omega \pi x_3) \Delta x^2  \\
		-\omega^2 \pi^2 \sin(\omega \pi x_4) \Delta x^2  \\
		-\omega^2 \pi^2 \sin(\omega \pi x_5) \Delta x^2  \\
		\omega \pi \cos(\omega \pi x_6) \Delta x  \\
	\end{bmatrix}
	\label{Eq:PPE 1D FDM matrix}
\end{equation}
For FDM, adding the individual equations gives zero on the left side but the right side need not be zero. Note that \cref{Eq:PPE 1D FDM matrix} uses a first order one sided stencil on the boundary to implement the Neumann condition. If a second order stencil is used, the left hand side of the equations also do not sum to zero.

\subsection{Detailed Analysis for Two Dimensional Problems}
Incompressible fluid flow equations are commonly solved by pressure projection methods which require solution of the pressure Poisson equation. Hence, in this section, we first analyze the Poisson equation in two dimensions with all Neumann boundary conditions:
\begin{equation}
\nabla^2 p(x,y) = \frac{\partial^2 p}{\partial x^2} + \frac{\partial^2 p}{\partial y^2} = S(x,y)
\label{Eq:Poisson eqn}
\end{equation}
In order to perform a rigorous error analysis, we use the method of manufactured solution. The source term $S(x,y)$ is computed for a given exact solution $p_{exact}(x,y)$:
\begin{equation}
\begin{aligned}
p_{exact}(x,y)&=\sin(\omega \pi x) \sin(\omega \pi y)\\
S(x,y) &= \nabla^2p_{exact} = -2\omega^2 \pi^2 \sin(\omega \pi x) \sin(\omega \pi y)
\end{aligned}
\label{Eq:Poisson eqn soln}
\end{equation}
where, $\omega$ is the wave number. The Poisson equation is solved over a unit square with Neumann boundary conditions. The normal derivative of the exact solution is imposed on all the four boundaries of the square.
\par For comparison, we first use the finite difference (FDM) and volume (FVM) methods on a unit square discretized with a structured grid. For this section, the wave number ($\omega$ in \cref{Eq:Poisson eqn soln}) is set to unity. FDM and FVM are the most popular methods used for computational fluid dynamics and heat transfer problems.
For FDM, we discretize the unit square using a structured Cartesian point placement. The Laplace operator in \cref{Eq:Poisson eqn} is approximated with the second order central difference stencil:
\begin{equation}
\begin{aligned}
\nabla^2 p(x,y) \approx &\frac{p(x+\Delta x, y) - 2 p(x,y) + p(x-\Delta x, y)}{\Delta x^2}\\
+ & \frac{p(x, y+\Delta y) - 2 p(x,y) + p(x, y-\Delta y)}{\Delta y^2}
\end{aligned}
\label{Eq:Poisson eqn FDM Laplace stencil}
\end{equation}
We use a uniform isotropic point placement with $\Delta x = \Delta y$. In order to impose the Neumann boundary condition on the four walls, we use a second order three-point stencil. For example, on the left wall, the first derivative can be approximated as:
\begin{equation}
\frac{\partial p}{\partial x} \approx \frac{-3 p(x,y) + 4 p(x+\Delta x, y) - p(x+2\Delta x, y)}{2 \Delta x}
\label{Eq:Poisson eqn FDM 1st derivative stencil}
\end{equation}

Now we briefly describe the finite volume method (FVM). Here, we discretize the unit square with uniform square shaped control volumes. In FVM, the governing \cref{Eq:Poisson eqn} is first integrated over a square control volume followed by application of the divergence theorem:
\begin{equation}
\sum_{f}\left( \nabla p \cdot \bm{N}\right) \big|_f = \sum_{f} \left( \nabla p_{exact} \cdot \bm{N} \right) \big|_f
\label{Eq:Poisson eqn FVM}
\end{equation}
where, $f$ denotes the face of the control volume and $\bm{N}$ is the unit normal of the face. The gradient at the face center is estimated using the central difference method based on the cell centered values of pressure. First order differencing is used on all the four walls to impose the Neumann boundary condition. This is a consistent formulation due to which all the discrete equations sum to zero. Thus, the column sum of the matrix and the sum of the sources for all the control volumes are zero when the Neumann condition is imposed at all boundaries.

\begin{figure}[H]
	\centering
	\begin{subfigure}[t]{0.45\textwidth}
		\includegraphics[width=\textwidth]{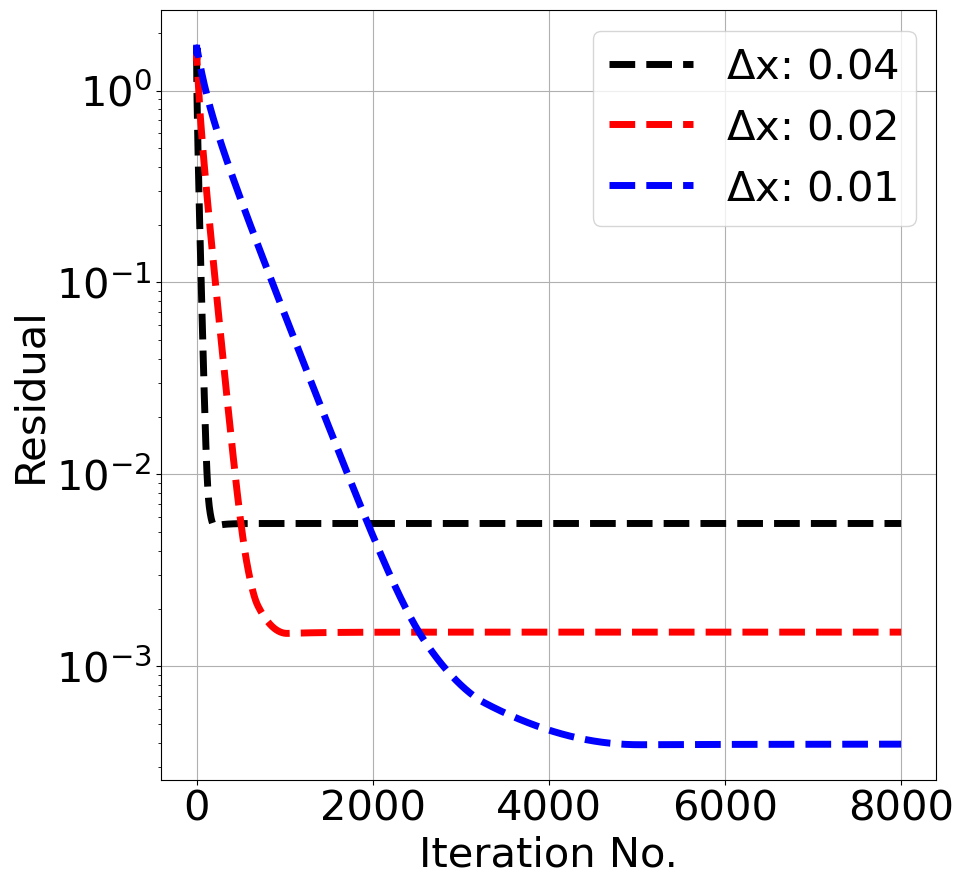}
		\caption{Finite Difference Method (FDM)}
		\label{Fig:Structured: res vs iter FDM}
	\end{subfigure}
	\begin{subfigure}[t]{0.45\textwidth}
		\includegraphics[width=\textwidth]{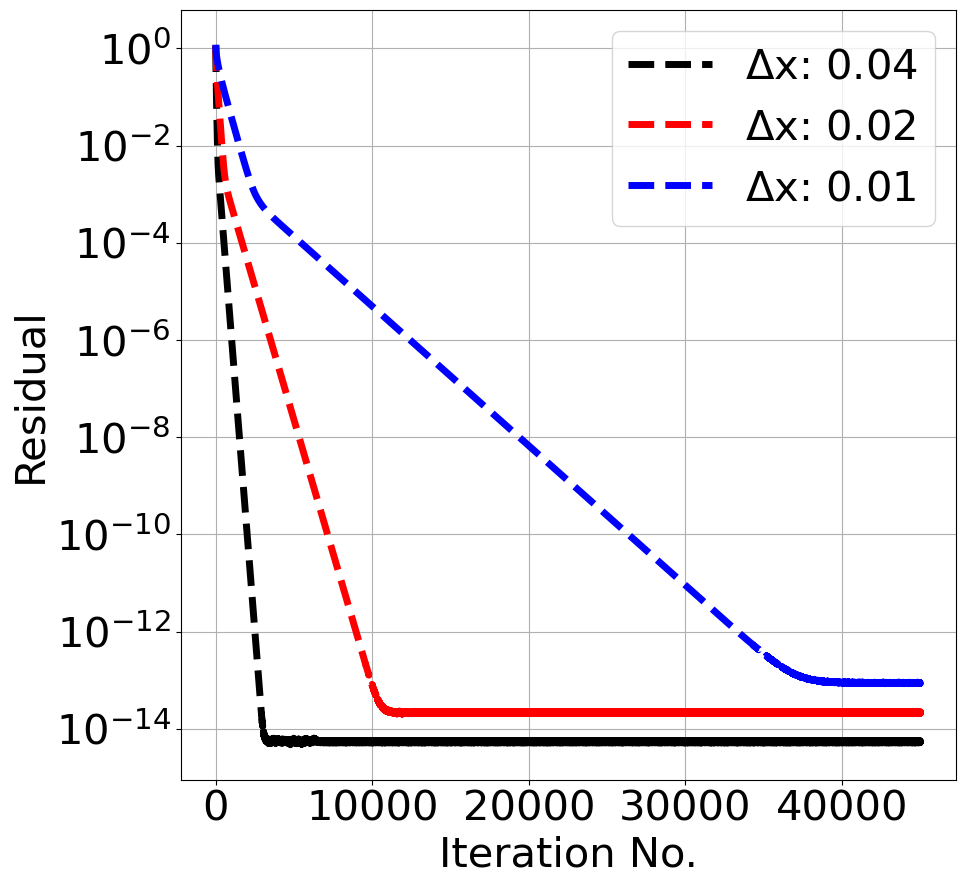}
		\caption{Finite Volume Method (FVM)}
		\label{Fig:Structured: res vs iter FVM}
	\end{subfigure}
	\caption{Residuals for Varying Grid Sizes}
	\label{Fig:Structured: res vs iter}
\end{figure}

\begin{figure}[H]
	\centering
	\begin{subfigure}[t]{0.45\textwidth}
		\includegraphics[width=\textwidth]{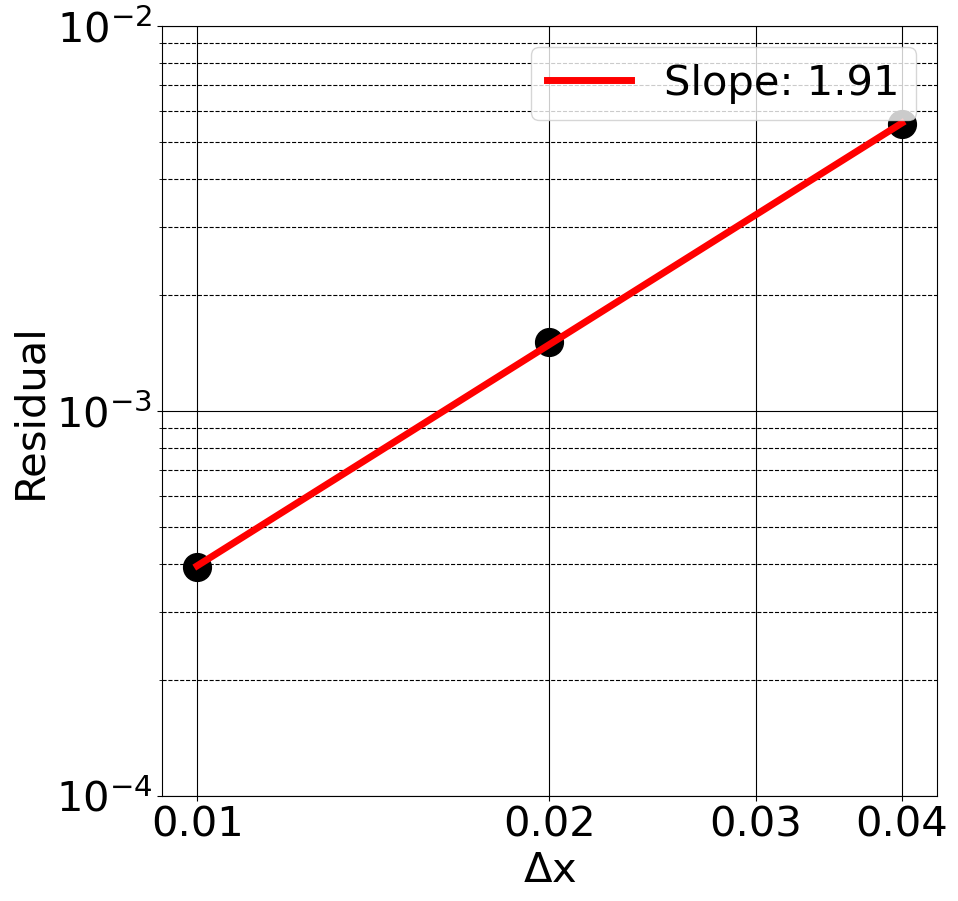}
		\caption{Final Residual for Finite Difference Method (FDM)}
		\label{Fig:Structured: rate of conv residual FDM}
	\end{subfigure}
	\begin{subfigure}[t]{0.45\textwidth}
		\includegraphics[width=\textwidth]{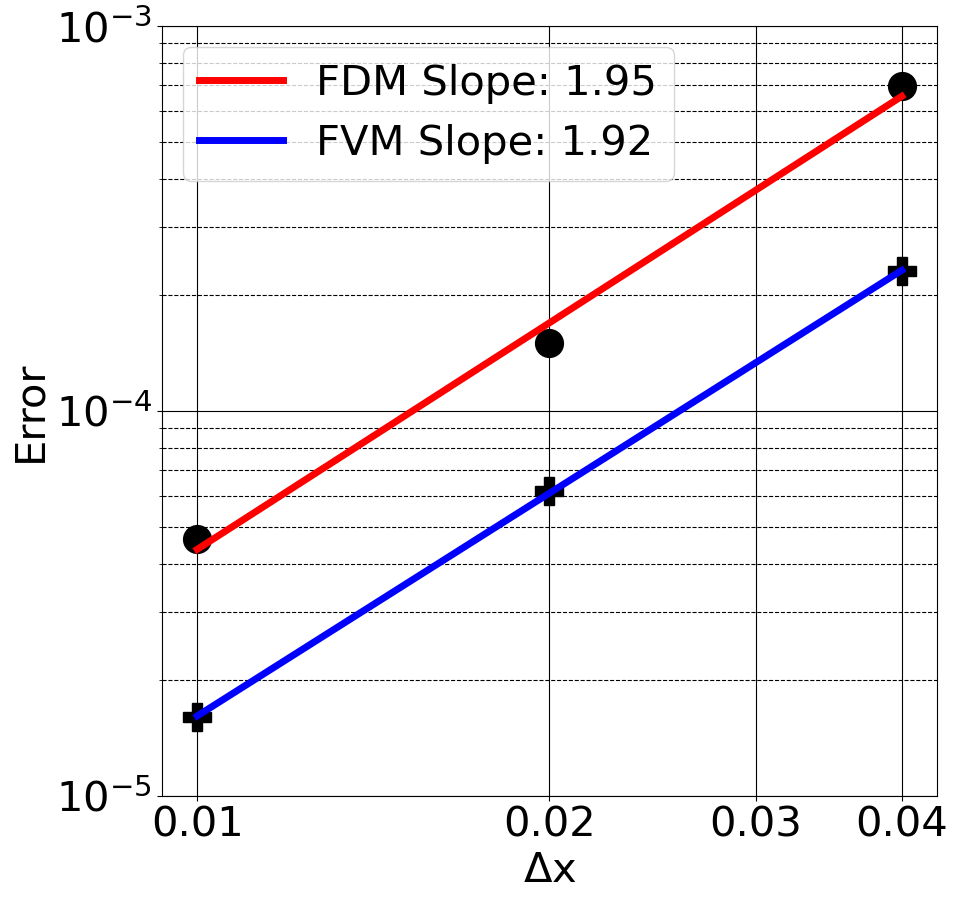}
		\caption{Solution Error for Finite Difference (FDM) and Volume (FVM) Methods}
		\label{Fig:Structured: rate of conv error}
	\end{subfigure}
	\caption{Rates of Convergence}
	\label{Fig:Structured: rate of conv}
\end{figure}

Both FDM and FVM discretizations lead to a linear equation $\bm{A}p=b$. The L$_1$ norms of relative residual and solution error are defined as follows:
\begin{equation}
\begin{aligned}
\text{Relative residual: }||r||_1 &= \frac{||\bm{A} p - b||_1}{||b||_1}\\
\text{Solution error: }||e||_1 &= ||p - p_{exact}||_1
\end{aligned}
\label{Eq:rel res error norms}
\end{equation}
where, $p_{exact}$ is given by \cref{Eq:Poisson eqn soln}. To analyze the residual and solution errors, three different point placements are used with $\Delta x = [0.04, 0.02, 0.01]$. This corresponds to [26, 51, 101] FDM points and [27, 52, 102] FVM control volumes in X and Y directions each. The FVM control volumes include the line segments at the boundaries used to impose the boundary conditions. The linear equation $\bm{A}p=b$ is solved by the successive over relaxation (SOR) method with an over-relaxation of 1.4. \Cref{Fig:Structured: res vs iter} plots the relative residuals against iteration number for both the methods. As expected, we observe that coarser point placement converges faster than a finer grid. The FDM does not explicitly impose the divergence theorem. The error in satisfaction of the divergence theorem is analogous to the truncation error which reduces with the grid refinement. Hence, the finest grid has lowest final residual. FVM is formulated in a flux conserving manner. This is because the flux that exits a control volume through a face enters the neighboring volume since the face cannot act as a sink or a source of that variable. Thus, the FVM inherently satisfies the divergence theorem. This distinction is reflected in the residual plots. The FDM residual saturates at a value based on grid refinement whereas, the FVM residual reduces to the round-off error.
Now we analyze the rate of convergence of the final saturated value of the FDM residual with grid refinement. \Cref{Fig:Structured: rate of conv residual FDM} plots the final value of the residual versus $\Delta x$ for three different refinements on a logarithmic scale. The slope of the best fit line (1.91) indicates that rate of convergence of the residual is second order. \Cref{Fig:Structured: rate of conv error} plots the convergence of the solution error for FDM and FVM. In this example, the actual value of the error is higher for the FDM than the FVM error. Both the methods show second order rates of convergence in solution errors. This is expected since we have used second order stencils for the discretization of the governing equation and the boundary condition (\cref{Eq:Poisson eqn FDM Laplace stencil,Eq:Poisson eqn FDM 1st derivative stencil,Eq:Poisson eqn FVM}). This exercise shows that unlike FVM, FDM satisfies divergence theorem only up to the accuracy of discretization error. Even though the residuals do not converge to round-off errors, the FDM solution error is still second order accurate. Hence, the FDM iterations can be stopped when the residual reaches a stationary value without increasing the solution error considerably. The solution error still displays the expected order of accuracy. We expect that if higher order finite difference discretizations are used consistently in the interior and at the boundary the same observations will be found with the stationary residuals and the solution errors decreasing per the order of discretization. This leads us to the study of the meshless method that strives to achieve high order convergence on complex domains. The rest of the paper describes the method and the observed rates of convergence in complex flows.

%

\section{Numerical Method} \label{Sec:Numerical Method}
In this section, we briefly describe the PHS-RBF approach. In the PHS-RBF method, we first discretize a complex domain with scattered points which do not have any connectivity in the form of a grid. Any scalar variable $s$ is then interpolated over a set of $q$ neighboring points known as a cloud using the formula,
\begin{equation}
s(\bm{x}) = \sum_{i=1}^{q} \lambda_i \phi_i (||\bm{x} - \bm{x_i}||_2) + \sum_{i=1}^{m} \gamma_i P_i (\bm{x})
\label{Eq:RBF_interp}
\end{equation}
where, the PHS-RBF is a function of radial distance ($r$) between the points: $\phi(r)=r^{2a+1}$ (here, we use $r^3$ function with $a=1$). $\lambda_i$ and $\gamma_i$ are $q+m$ unknown coefficients and $P(\bm{x})$ is a set of appended monomials. For a $d$ dimensional problem (which may be 1, 2 or 3) and maximum degree of the appended monomials as $k$, the number of appended monomials is given by the binomial coefficient $m=\binom{k+d}{k}$. For example, if $d=2$ and $k=2$, there are $m=\binom{k+d}{k}=\binom{2+2}{2}=6$ monomial terms given by $[1, x, y, x^2, xy, y^2]$. \Cref{Eq:RBF_interp} is collocated at the $q$ cloud points to give $q$ equations. Since there are a total of $q+m$ unknown coefficients, additional $m$ constraints are required to close the system of equations \cite{flyer2016onrole_I}. These are given by
\begin{equation}
\sum_{i=1}^{q} \lambda_i P_j(\bm{x_i}) =0 \hspace{0.5cm} \text{for } 1 \leq j \leq m
\label{Eq:RBF_constraint}
\end{equation}
Writing \cref{Eq:RBF_interp,Eq:RBF_constraint} in a matrix vector form gives:
\begin{equation}
\begin{bmatrix}
\bm{\Phi} & \bm{P}  \\
\bm{P}^T & \bm{0} \\
\end{bmatrix}
\begin{bmatrix}
\bm{\lambda}  \\
\bm{\gamma} \\
\end{bmatrix} =
\begin{bmatrix}
\bm{A}
\end{bmatrix}
\begin{bmatrix}
\bm{\lambda}  \\
\bm{\gamma} \\
\end{bmatrix} =
\begin{bmatrix}
\bm{s}  \\
\bm{0} \\
\end{bmatrix}
\label{Eq:RBF_interp_mat_vec}
\end{equation}
where, the superscript $T$ denotes transpose, $\bm{\lambda} = [\lambda_1,...,\lambda_q]^T$, $\bm{\gamma} = [\gamma_1,...,\gamma_m]^T$, $\bm{s} = [s(\bm{x_1}),...,s(\bm{x_q})]^T$ and $\bm{0}$ is the matrix of all zeros. Sizes of the submatrices $\bm{\Phi}$ and $\bm{P}$ are $q\times q$ and $q\times m$ respectively. Evaluating the PHS-RBF and monomials at the cloud points gives these submatrices:
\begin{equation}
\bm{\Phi} =
\begin{bmatrix}
\phi \left(||\bm{x_1} - \bm{x_1}||_2\right) & \dots  & \phi \left(||\bm{x_1} - \bm{x_q}||_2\right) \\
\vdots & \ddots & \vdots \\
\phi \left(||\bm{x_q} - \bm{x_1}||_2\right) & \dots  & \phi \left(||\bm{x_q} - \bm{x_q}||_2\right) \\
\end{bmatrix}
\label{Eq:RBF_interp_phi}
\end{equation}
\begin{equation}
\bm{P} =
\begin{bmatrix}
1 & x_1  & y_1 & x_1^2 & x_1 y_1 & y_1^2 \\
\vdots & \vdots & \vdots & \vdots & \vdots & \vdots \\
1 & x_q  & y_q & x_q^2 & x_q y_q & y_q^2 \\
\end{bmatrix}
\label{Eq:RBF_interp_poly}
\end{equation}
In order to interpolate any function $s$, \cref{Eq:RBF_interp_mat_vec} is solved for the $q+m$ unknown coefficients $\lambda_i$ and $\gamma_i$.
\begin{equation}
\begin{bmatrix}
\bm{\lambda}  \\
\bm{\gamma} \\
\end{bmatrix} =
\begin{bmatrix}
\bm{A}
\end{bmatrix} ^{-1}
\begin{bmatrix}
\bm{s}  \\
\bm{0} \\
\end{bmatrix}
\label{Eq:RBF_interp_mat_vec_solve}
\end{equation}
where, inverse of the matrix $\bm{A}$ is shown only as a notation in this derivation. In actual computations, we use a linear solver.
\par The Navier-Stokes and Poisson equations have differential operators such as gradient and Laplacian. Hence, it is important to compute these derivatives. These can be estimated as a weighted sum of the function values evaluated at the cloud points. Let $\mathcal{L}$ denote a scalar linear operator such as first or second derivative. Operating $\mathcal{L}$ on \cref{Eq:RBF_interp} and using the linearity gives:

\begin{equation}
\mathcal{L} [s(\textbf{x})] = \sum_{i=1}^{q} \lambda_i \mathcal{L} [\phi_i (\bm{x})] + \sum_{i=1}^{m} \gamma_i \mathcal{L}[P_i (\bm{x})]
\label{Eq:RBF_interp_L}
\end{equation}
\Cref{Eq:RBF_interp_L} is evaluated at the $q$ cloud points giving a rectangular matrix vector system:
\begin{equation}
\mathcal{L}[\bm{s}] =
\begin{bmatrix}
\mathcal{L}[\bm{\Phi}] & \mathcal{L}[\bm{P}]  \\
\end{bmatrix}
\begin{bmatrix}
\bm{\lambda}  \\
\bm{\gamma} \\
\end{bmatrix}
\label{Eq:RBF_interp_mat_vec_L}
\end{equation}
where, $\mathcal{L}[\bm{\Phi}]$ and $\mathcal{L}[\bm{P}]$ are matrices of sizes $q\times q$ and $q\times m$ respectively.
\begin{equation}
\mathcal{L}[\bm{s}] = [\mathcal{L}[s(\bm{x})]_{\bm{x_1}},...,\mathcal{L}[s(\bm{x})]_{\bm{x_q}}]^T
\label{Eq:RBF_interp_phi_Ls}
\end{equation}
\begin{equation}
\mathcal{L}[\bm{\Phi}] =
\begin{bmatrix}
\mathcal{L}[\phi \left(||\bm{x} - \bm{x_1}||_2\right)]_{\bm{x_1}} & \dots  & \mathcal{L}[\phi \left(||\bm{x} - \bm{x_q}||_2\right)]_{\bm{x_1}} \\
\vdots & \ddots & \vdots \\
\mathcal{L}[\phi \left(||\bm{x} - \bm{x_1}||_2\right)]_{\bm{x_q}} & \dots  & \mathcal{L}[\phi \left(||\bm{x} - \bm{x_q}||_2\right)]_{\bm{x_q}} \\
\end{bmatrix}
\label{Eq:RBF_interp_phi_L}
\end{equation}
Similar to \cref{Eq:RBF_interp_poly}, for $d=2$ and $k=2$, $\mathcal{L}[\bm{P}]$ is given by:
\begin{equation}
\mathcal{L}[\bm{P}] =
\begin{bmatrix}
\mathcal{L}[1]_{\bm{x_1}} & \mathcal{L}[x]_{\bm{x_1}}  & \mathcal{L}[y]_{\bm{x_1}} & \mathcal{L}[x^2]_{\bm{x_1}} & \mathcal{L}[x y]_{\bm{x_1}} & \mathcal{L}[y^2]_{\bm{x_1}} \\
\vdots & \vdots & \vdots & \vdots & \vdots & \vdots \\
\mathcal{L}[1]_{\bm{x_q}} & \mathcal{L}[x]_{\bm{x_q}}  & \mathcal{L}[y]_{\bm{x_q}} & \mathcal{L}[x^2]_{\bm{x_q}} & \mathcal{L}[x y]_{\bm{x_q}} & \mathcal{L}[y^2]_{\bm{x_q}} \\
\end{bmatrix}
\label{Eq:RBF_interp_poly_L}
\end{equation}
The subscripts in \cref{Eq:RBF_interp_phi_L,Eq:RBF_interp_poly_L,Eq:RBF_interp_phi_Ls} indicate that the matrices are obtained by operating $\mathcal{L}$ on the basis functions and evaluating them at the cloud points ($\bm{x_i}$). Substituting the values of coefficient from \cref{Eq:RBF_interp_mat_vec_solve} in \cref{Eq:RBF_interp_mat_vec_L} and simplifying, we get:
\begin{equation}
\begin{aligned}
\mathcal{L}[\bm{s}] &=
\left(\begin{bmatrix}
\mathcal{L}[\bm{\Phi}] & \mathcal{L}[\bm{P}]  \\
\end{bmatrix}
\begin{bmatrix}
\bm{A}
\end{bmatrix} ^{-1}\right)
\begin{bmatrix}
\bm{s}  \\
\bm{0} \\
\end{bmatrix}
=
\begin{bmatrix}
\bm{B}
\end{bmatrix}
\begin{bmatrix}
\bm{s}  \\
\bm{0} \\
\end{bmatrix}\\
&=
\begin{bmatrix}
\bm{B_1} & \bm{B_2}
\end{bmatrix}
\begin{bmatrix}
\bm{s}  \\
\bm{0} \\
\end{bmatrix}
= [\bm{B_1}] [\bm{s}] + [\bm{B_2}] [\bm{0}]
= [\bm{B_1}] [\bm{s}]
\end{aligned}
\label{Eq:RBF_interp_mat_vec_L_solve}
\end{equation}
Splitting the matrix $[\bm{B}]$ of size $q \times (q+m)$ gives two submatrices $[\bm{B_1}]$ and $[\bm{B_2}]$ of sizes $q \times q$ and $q \times m$ respectively. Matrix $[\bm{B_1}]$ contains the coefficients used in linear combination of the function at the cloud points to estimate the operator ($\mathcal{L}$) values at those points. These coefficients at each point in the entire domain are assembled into a sparse matrix. Multiplying the sparse matrix by the velocity or pressure in the domain gives a numerical estimate of the gradient or Laplacian. By satisfying the governing equation at all the scattered points gives a set of discrete equations for the unknown values of the variable at the scattered points. For further details about the PHS-RBF coefficients and their properties such as condition number and accuracy, please refer to our previous work \cite{shahane2021high}.

\section{Poisson Equation with Manufactured Solution using Meshless Discretization} \label{Sec:Poisson Equation with Manufactured Solution Meshless Discretization}

\begin{figure}[H]
	\centering
	\begin{subfigure}[t]{0.45\textwidth}
		\includegraphics[width=\textwidth]{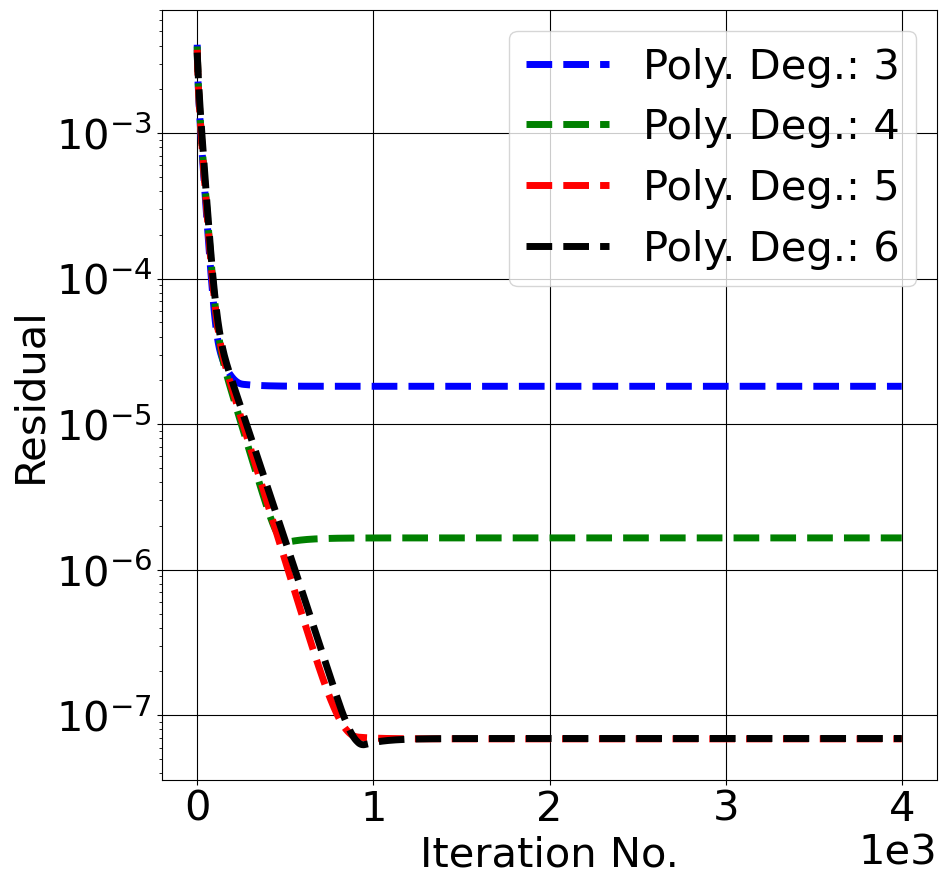}
		\caption{Average $\Delta x=0.0418$}
		\label{Fig:Meshless: poisson res vs iter wv 1 coarse}
	\end{subfigure}
	\hspace{0.05\textwidth}
	\begin{subfigure}[t]{0.45\textwidth}
		\includegraphics[width=\textwidth]{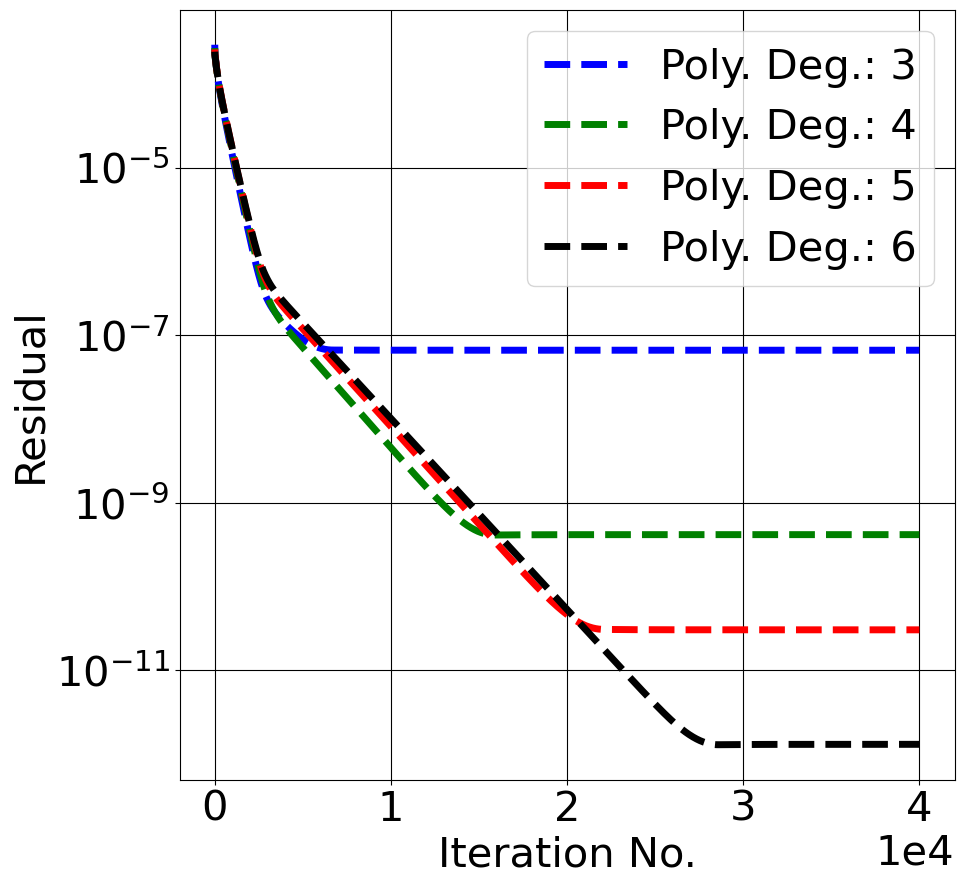}
		\caption{Average $\Delta x=0.0101$}
		\label{Fig:Meshless: poisson res vs iter wv 1 fine}
	\end{subfigure}
	\caption{Residuals for Varying Grid Sizes for Wave Number $\omega=1$}
	\label{Fig:Meshless: poisson res vs iter wv 1}
\end{figure}

\begin{figure}[H]
	\centering
	\begin{subfigure}[t]{0.45\textwidth}
		\includegraphics[width=\textwidth]{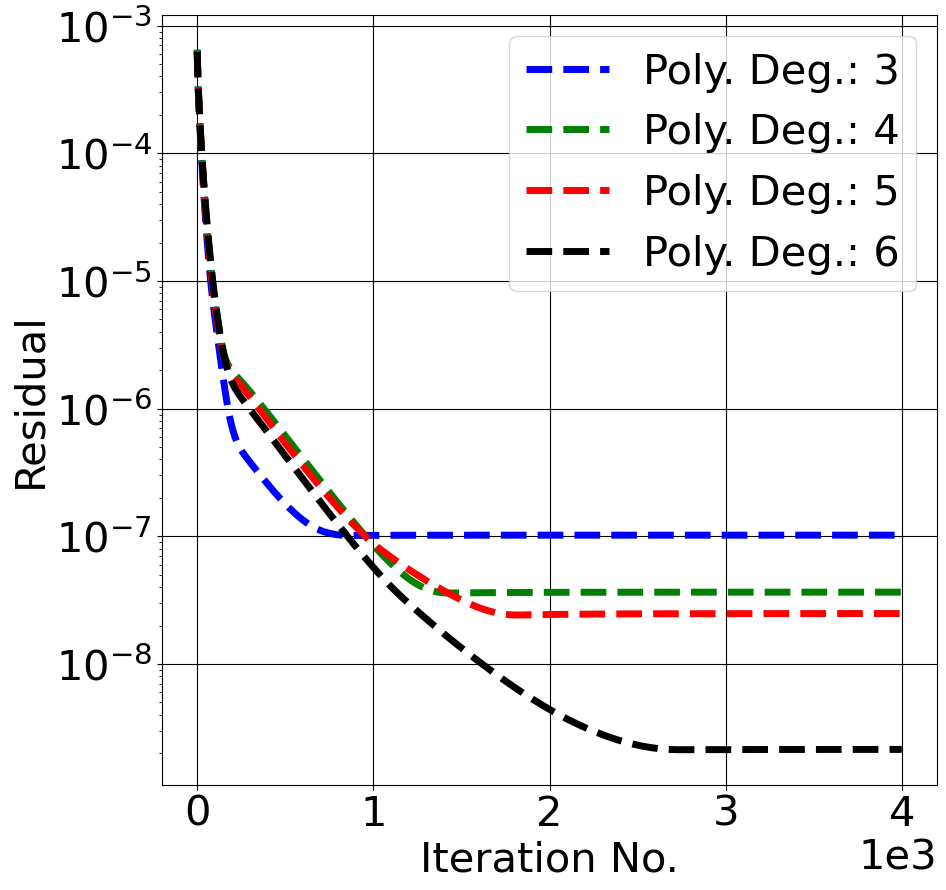}
		\caption{Average $\Delta x=0.0203$}
		\label{Fig:Meshless: poisson res vs iter wv 4 coarse}
	\end{subfigure}
	\hspace{0.05\textwidth}
	\begin{subfigure}[t]{0.45\textwidth}
		\includegraphics[width=\textwidth]{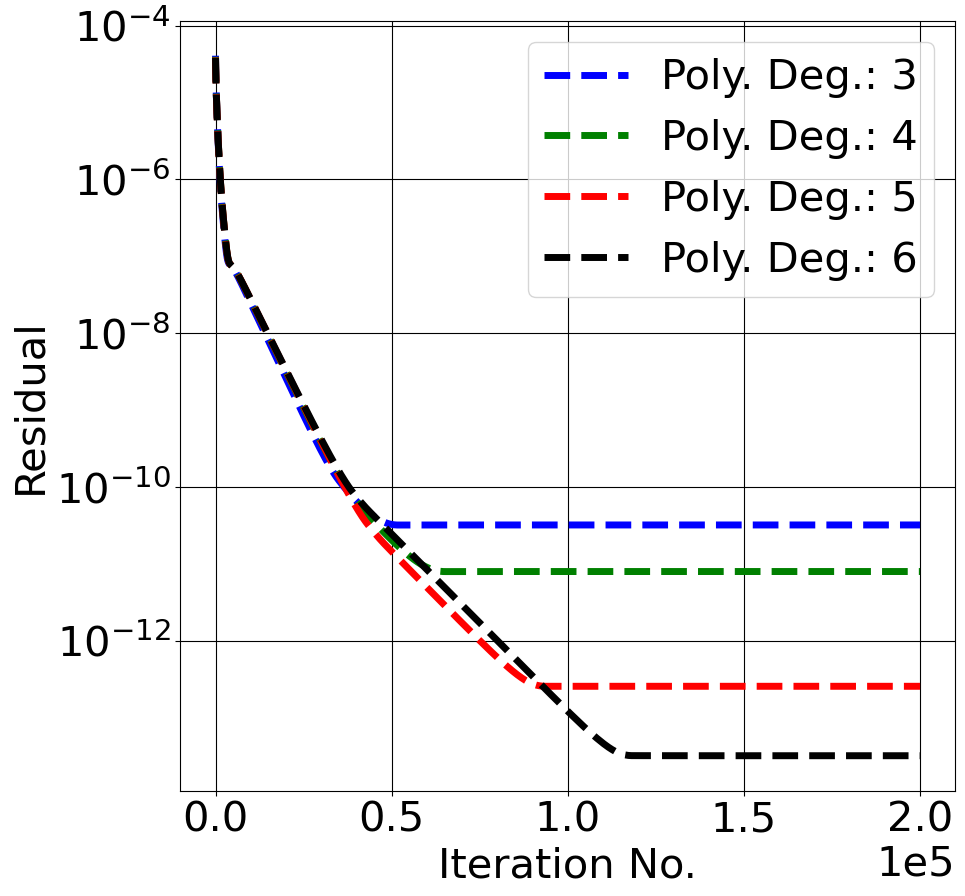}
		\caption{Average $\Delta x=0.0050$}
		\label{Fig:Meshless: poisson res vs iter wv 4 fine}
	\end{subfigure}
	\caption{Residuals for Varying Grid Sizes for Wave Number $\omega=4$}
	\label{Fig:Meshless: poisson res vs iter wv 4}
\end{figure}

The higher order meshless method discussed in \cref{Sec:Numerical Method} is now applied to solve the Poisson \cref{Eq:Poisson eqn} with all Neumann boundary conditions for manufactured solution with two values of wave numbers: $\omega=1$ and $\omega=4$ (\cref{Eq:Poisson eqn soln}). The PHS-RBF stencils of the Laplacian and gradient operators are used for the interior and boundary points respectively. The unit square is discretized with an unstructured triangular grid using the open source software GMSH \cite{geuzaine2009gmsh} and the vertices of the triangular elements are then used as scattered points. For the cases of $\omega=1$ and $\omega=4$, three point placements with an average $\Delta x$ of $[0.0418, 0.0203, 0.0101]$ and $[0.0203, 0.0101, 0.0050]$ respectively are used. Degrees of the appended polynomials with the PHS-RBF are varied from 3 to 6 in order to increase the order of accuracy. The discrete linear equation of the form $\bm{A}p=b$ is solved iteratively with the SOR algorithm with an over-relaxation of 1.4. \Cref{Fig:Meshless: poisson res vs iter wv 1,Fig:Meshless: poisson res vs iter wv 4} plot the relative residuals (\cref{Eq:rel res error norms}) versus SOR iterations for the coarsest and finest grid of both the wave numbers. Since the meshless method is akin to a high order accurate FDM, we observe that in this case also the residual saturates at a value much higher than the round-off error. The case of higher wave number plotted in \cref{Fig:Meshless: poisson res vs iter wv 4} requires more iterations than the low wave number case (\cref{Fig:Meshless: poisson res vs iter wv 1}) since more number of points are used to get higher resolution \cite{shahane2021high}. Further, as we increase the degree of the polynomial, the order of discretization convergence and the order of residual convergence increase. Thus, the saturation value of the residual decreases with the polynomial degree. As we increase the degree of appended polynomial and the number of the scattered points (both h and p refinements), we see that the inconsistency eventually vanishes and does not affect the solution obtained. We also observe that using higher degree of polynomials for a fixed number of scattered points does not adversely affect the convergence of the SOR iterative solver. This is very encouraging despite the increase in the bandwidth of the coefficient matrix.

\begin{figure}[H]
	\centering
	\begin{subfigure}[t]{0.45\textwidth}
		\includegraphics[width=\textwidth]{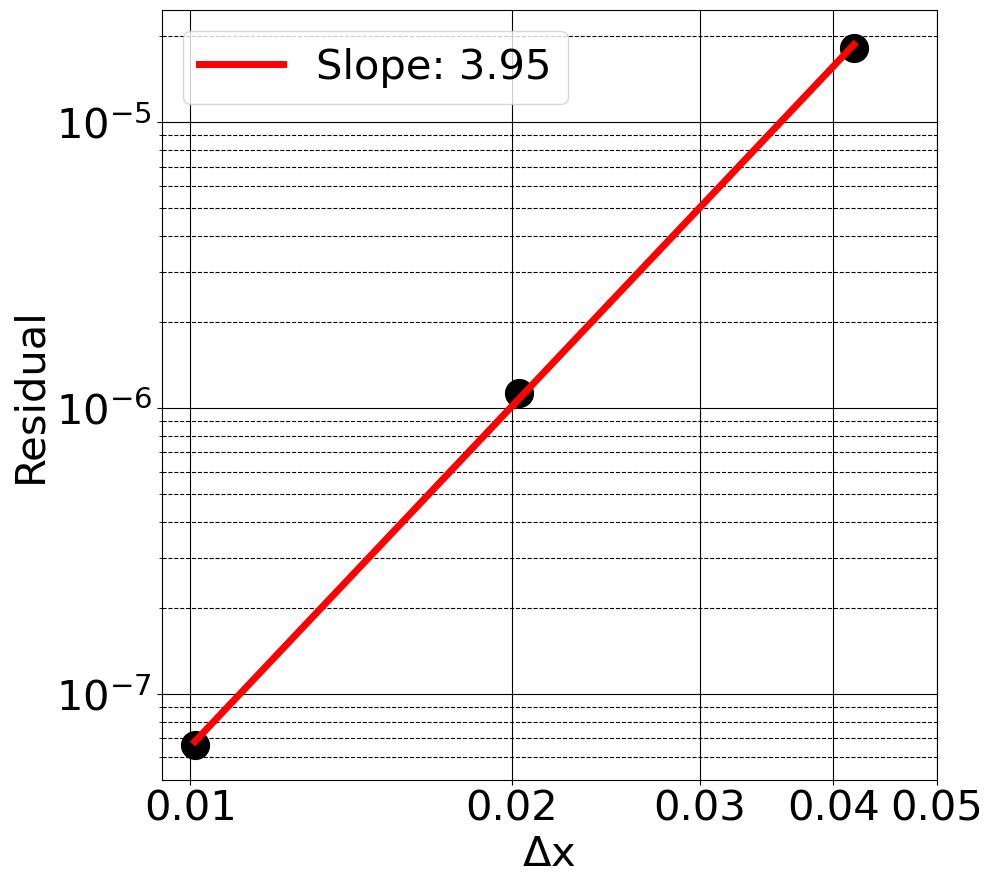}
		\caption{Degree of Appended Polynomial: 3}
	\end{subfigure}
	\hspace{0.05\textwidth}
	\begin{subfigure}[t]{0.45\textwidth}
		\includegraphics[width=\textwidth]{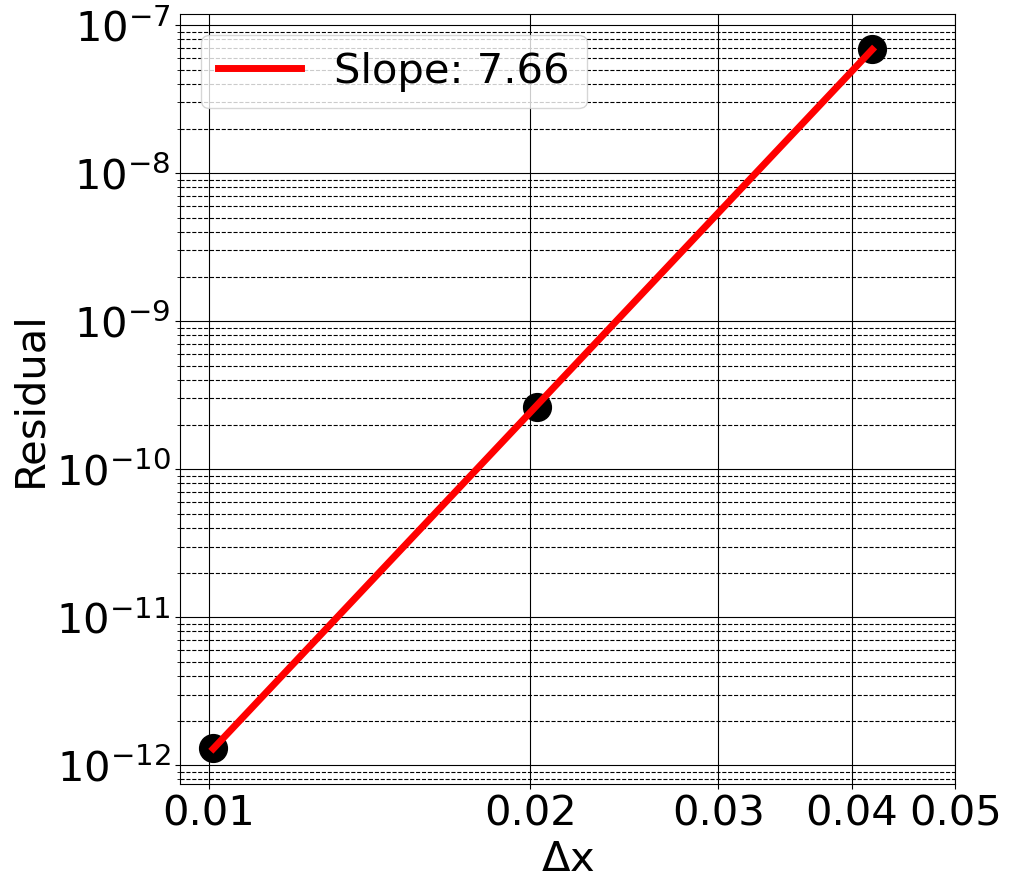}
		\caption{Degree of Appended Polynomial: 6}
	\end{subfigure}
	\caption{Convergence Rate of Final Residuals: Wave Number $\omega=1$}
	\label{Fig:Meshless: poisson rate of conv residual wv 1}
\end{figure}
\begin{figure}[H]
	\centering
	\begin{subfigure}[t]{0.45\textwidth}
		\includegraphics[width=\textwidth]{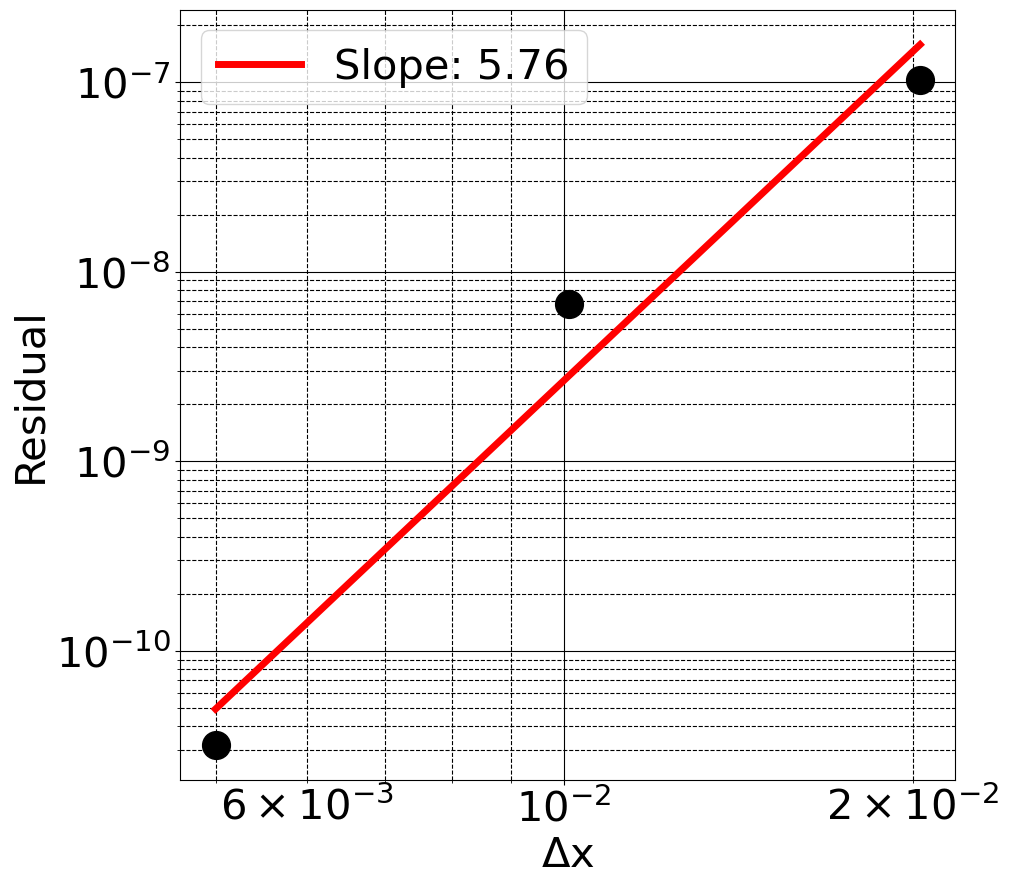}
		\caption{Degree of Appended Polynomial: 3}
	\end{subfigure}
	\hspace{0.05\textwidth}
	\begin{subfigure}[t]{0.45\textwidth}
		\includegraphics[width=\textwidth]{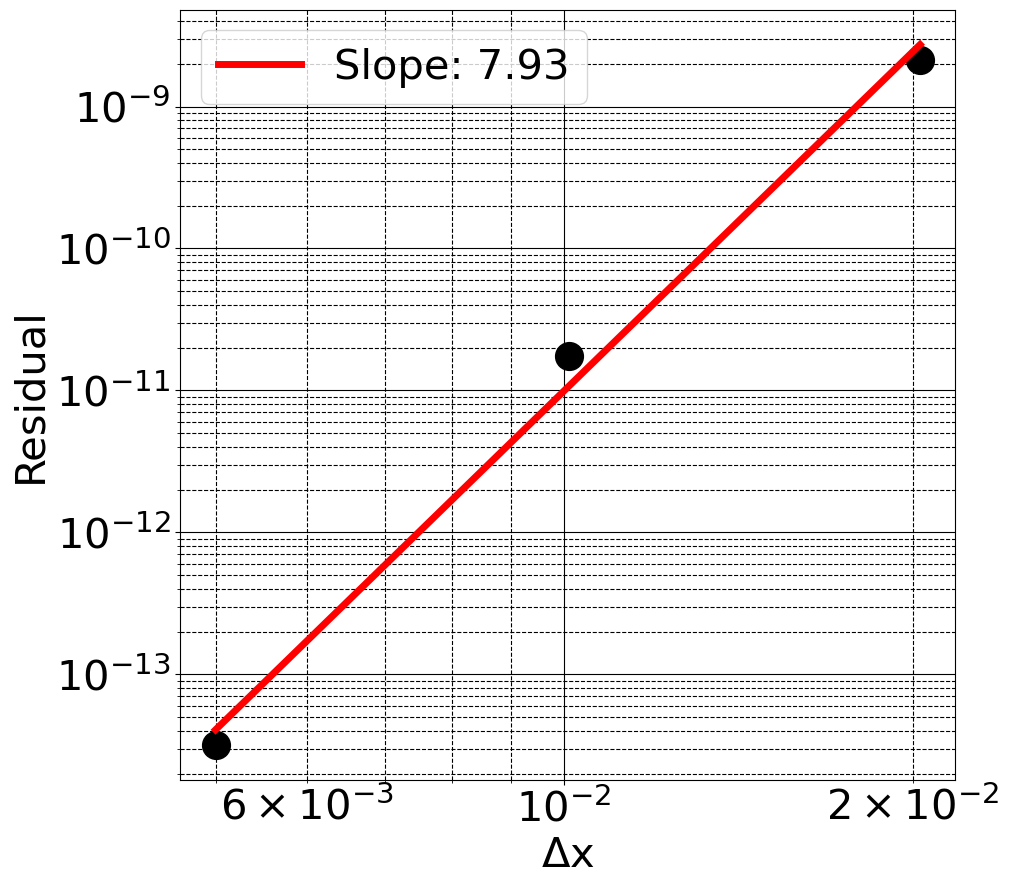}
		\caption{Degree of Appended Polynomial: 6}
	\end{subfigure}
	\caption{Convergence Rate of Final Residuals: Wave Number $\omega=4$}
	\label{Fig:Meshless: poisson rate of conv residual wv 4}
\end{figure}
\begin{figure}[H]
	\centering
	\includegraphics[width=0.6\textwidth]{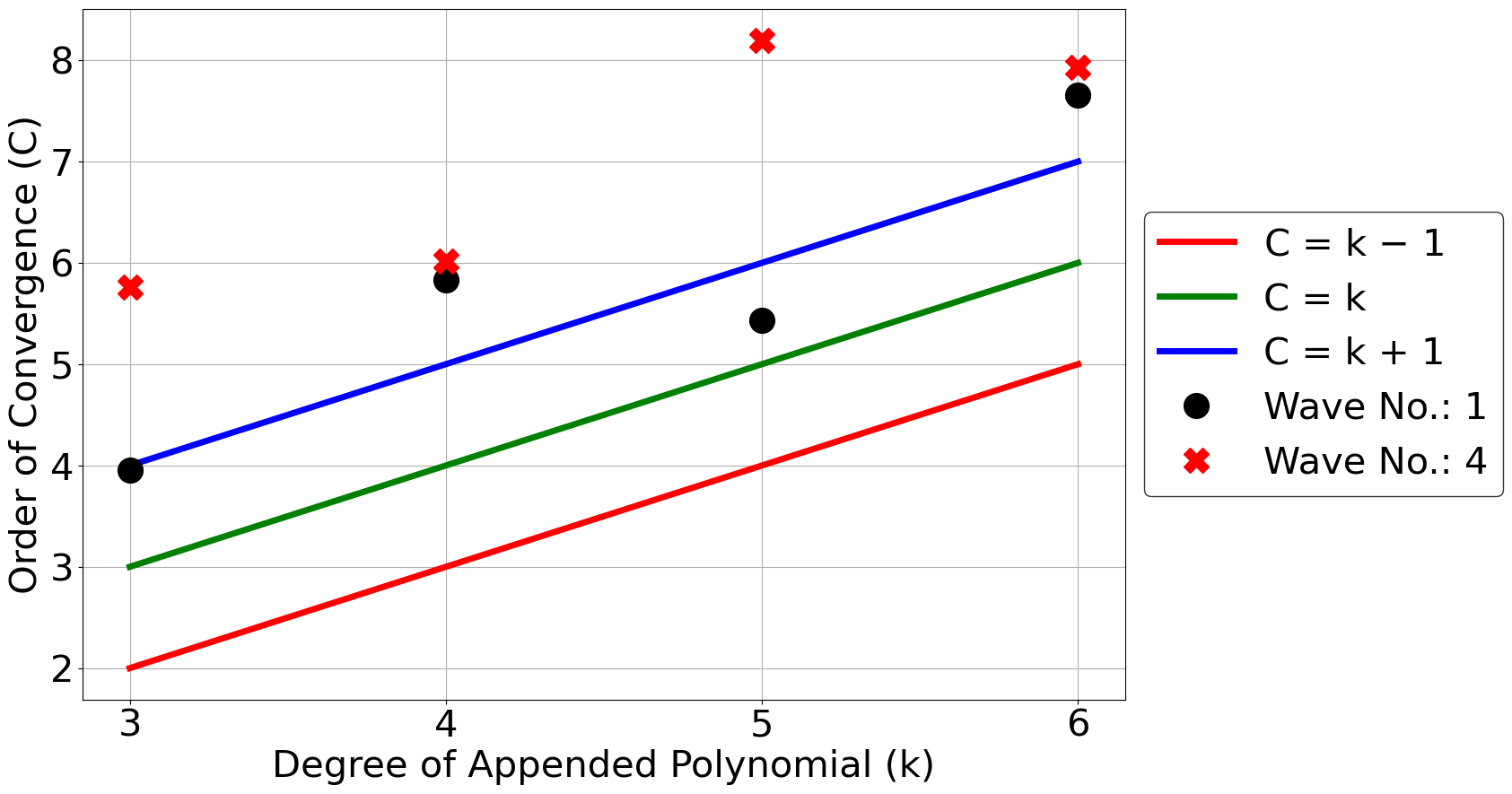}
	\caption{Composite Convergence Rate of Final Residuals with Degrees of Appended Polynomial}
	\label{Fig:Meshless: poisson rate of conv residual with polydeg}
\end{figure}

From the various simulations, we have again extracted the rates of convergence of the stationary residuals with average spacing of the points, for different degrees of polynomial. \Cref{Fig:Meshless: poisson rate of conv residual wv 1,Fig:Meshless: poisson rate of conv residual wv 4} plot the saturated residuals with average $\Delta x$ on a logarithmic scale. The best fit lines estimate the orders of convergence. For the sake of brevity, plots for only two polynomial degrees are included here. In \cref{Fig:Meshless: poisson rate of conv residual with polydeg}, we have combined all the slopes as a composited plot for both the wave numbers and four values of polynomial degrees. Three guide lines corresponding to $k-1$, $k$ and $k+1$ are also plotted for reference where, $k$ is the polynomial degree. It is expected that the discretization errors for a second derivative converge at the rate of $k-1$ or higher. Here, we observe that the saturation residuals converge much faster, typically $k+1$ or higher. This is a very encouraging and powerful result, as the inconsistency is seen to vanish much faster than the discretization errors if higher order polynomials are used. Since the point placement is irregular, the estimate of $\Delta x$ is approximate. Hence, we observe some oscillations in these rates. Moreover, all the rates are at least $k^\text{th}$ order accurate.

\begin{figure}[H]
	\centering
	\includegraphics[width=0.6\textwidth]{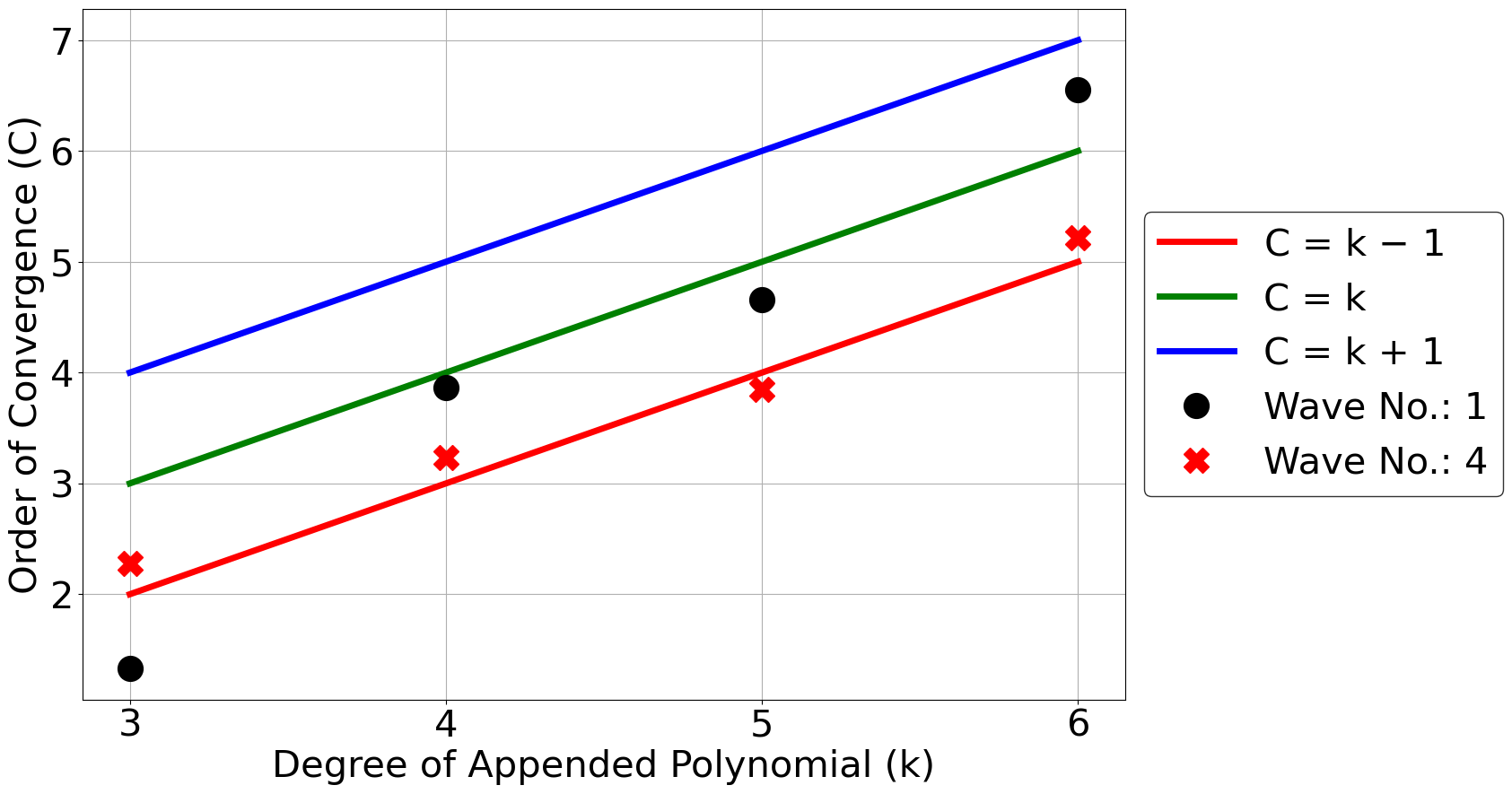}
	\caption{Composite Convergence Rate of Solution Error with Degrees of Appended Polynomial}
	\label{Fig:Meshless: poisson rate of conv error with polydeg}
\end{figure}

We next estimate the rate of convergence of the solution error (\cref{Eq:rel res error norms}) using the manufactured solution (\cref{Eq:Poisson eqn soln}). \Cref{Fig:Meshless: poisson rate of conv error with polydeg} shows that these rates also improve monotonically with the polynomial degrees as expected, for both the wave numbers. We have previously observed that the PHS-RBF stencils with appended polynomials of degree $k$ give $\mathcal{O}(k+1)$, $\mathcal{O}(k)$ and $\mathcal{O}(k-1)$ convergence for the interpolation, gradient and Laplacian operators respectively \cite{shahane2021high,flyer2016onrole_I}. Since the governing equation involves the Laplacian operator and the Neumann boundary has gradient operator, we observe that the rates lie between the guidelines of $k$ and $k-1$.

\begin{figure}[H]
	\centering
	\begin{subfigure}[t]{0.45\textwidth}
		\includegraphics[width=\textwidth]{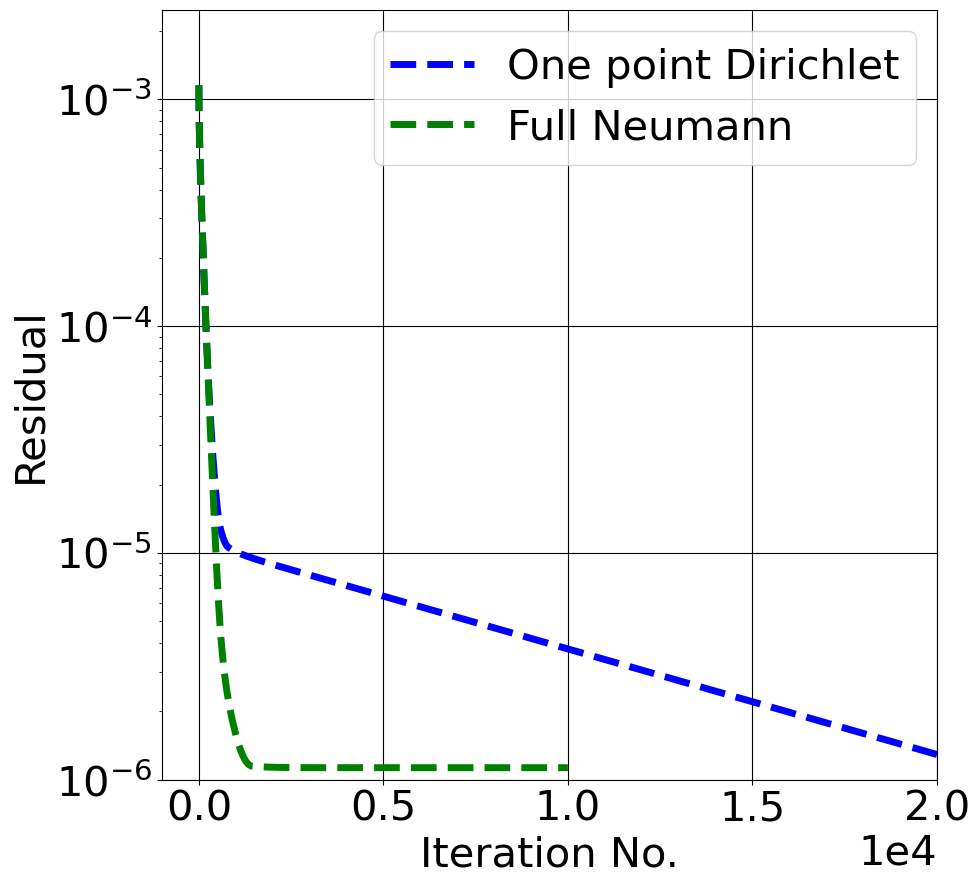}
		\caption{Degree of Appended Polynomial: 3}
	\end{subfigure}
	\hspace{0.05\textwidth}
	\begin{subfigure}[t]{0.45\textwidth}
		\includegraphics[width=\textwidth]{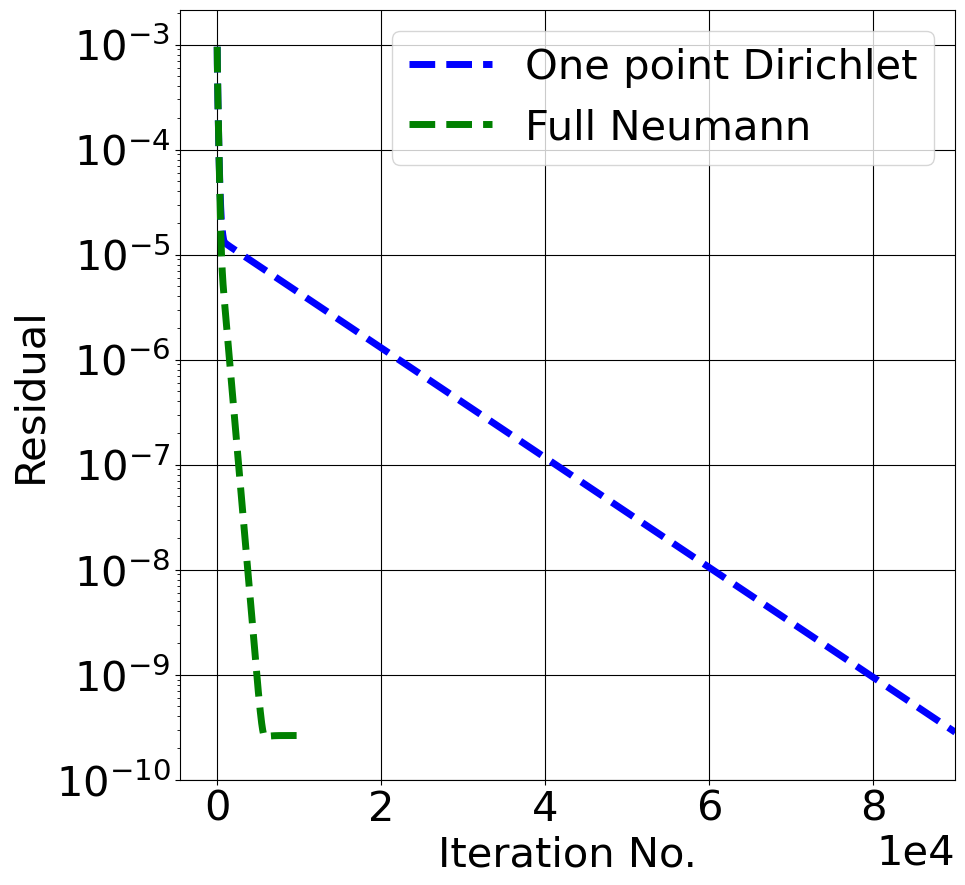}
		\caption{Degree of Appended Polynomial: 6}
	\end{subfigure}
	\caption{Comparison of Full Neumann and One Point Dirichlet Boundary Conditions (Wave Number $\omega=1$, Average $\Delta x=0.0203$)}
	\label{Fig:Meshless: poisson residual full_neu and 1_pt_dir}
\end{figure}
The Poisson equation with full Neumann boundary condition gives an ill-conditioned linear system (rank deficient matrix). The ill-conditioning can be resolved by prescribing a constraint on either the mean of the pressures or the pressure at a point. The former called the regularization method \cite{medusa_poisson_wiki,shahane2021high} fixes the mean to a constant value (zero in default) by adding an extra equation at the end of the coefficient matrix. This last row has a zero on the diagonal and presents difficulty in solution by point methods such as SOR. Methods such as BiCGSTAB internally rearrange the equations and are able to converge. Another approach is to fix the pressure at an arbitrary point to a reference value (as mentioned earlier). This also regularizes the equations but converges much slower than when the level is allowed to float. Relaxation schemes such as SOR can be used for such a constrained system with a fixed point pressure. We have tried prescribing a single boundary point with Dirichlet condition with all other points having Neumann condition. \Cref{Fig:Meshless: poisson residual full_neu and 1_pt_dir} compares the convergence of one point Dirichlet with full Neumann condition. We observe that the single point Dirichlet condition significantly slows the convergence. It requires five and ten times more iterations for polynomial degrees of 3 and 6 respectively to reach the stationary level of residuals of the full Neumann case. This shows that it is advantageous to iterate the original discrete equations with the full Neumann boundary condition than the regularized version of the original coefficient matrix. Further, we observe that iterating with the SOR solver until the residual saturates gives higher order convergence in the solution error similar to the BiCGSTAB solver with regularization \cite{shahane2021high}. Thus, we believe that the inconsistency and the ill-conditioning can be simply and efficiently resolved by iterating with an inexpensive solver such as SOR and terminating the iterations when the residual reaches the stationary level. The SOR solver can be easily coupled with a multilevel method to accelerate convergence. We are in the process of developing such an algorithm \cite{radhakrishnan2021non}

\section{Incompressible Fluid Flow Problems} \label{Sec:Incompressible Fluid Flow Problems}
We next study the convergence of a fractional step algorithm for incompressible flows in which the pressure Poisson equation is solved with a SOR \cite{shahane2021high}. We have considered four model problems and examine the level of residuals at stationarity and its convergence with point refinement and polynomial degree. Steady state solutions are obtained by time marching. We define the steady state convergence as  $\frac{\phi^{n+1} - \phi^n}{\Delta t}$ where, $n$ and $n+1$ are consecutive timesteps and $\phi$ denotes components of the velocity field. Steady state is assumed when the above error reaches a tolerance of 1E--7. Starting with a timestep $n$, intermediate velocity fields $\hat{u}$ and $\hat{v}$ are estimated by explicit Euler method:
\begin{equation}
\rho \frac{\hat{u} - u^n}{\Delta t} = -\rho \bm{u}^n \bullet (\nabla u^n) + \mu \nabla^2 u^n
\label{Eq:frac step u hat}
\end{equation}
\begin{equation}
\rho \frac{\hat{v} - v^n}{\Delta t} = -\rho \bm{u}^n \bullet (\nabla v^n) + \mu \nabla^2 v^n
\label{Eq:frac step v hat}
\end{equation}
where, $\rho$ is density, $\mu$ is viscosity and $\Delta t$ is timestep. Note that the above equations do not have the pressure gradient term. Subtracting \cref{Eq:frac step u hat,Eq:frac step v hat} from full momentum equations, applying the divergence operator and imposing the continuity equation gives the pressure Poisson equation:
\begin{equation}
\nabla \bullet (\nabla p) = \frac{\rho}{\Delta t} \left(\frac{\partial \hat{u}}{\partial x} + \frac{\partial \hat{v}}{\partial y} \right)
\label{Eq:frac step PPE}
\end{equation}
The above equation is solved with full Neumann boundary conditions consistent with boundary conditions of prescribed velocities. The pressure boundary condition is obtained by imposing the normal momentum equation:
\begin{equation}
\nabla p \bullet \bm{N} = (-\rho(\bm{u}\bullet \nabla)\bm{u} + \mu \nabla^2 \bm{u} )\bullet \bm{N}
\label{Eq:frac step normal momentum}
\end{equation}
which is derived by taking the dot product of the momentum equations with the unit normal $\bm{N}$ at the domain boundaries. The gradient and Laplacian operators in \cref{Eq:frac step PPE,Eq:frac step normal momentum,Eq:frac step u hat,Eq:frac step v hat} are estimated using the PHS-RBF stencils. Since momentum \cref{Eq:frac step u hat,Eq:frac step v hat} are explicit, a linear solver is required only for the pressure Poisson equation. Here, we use the SOR solver with the over-relaxation of 1.2. Since we are interested only in the steady state solutions, we use only 10-20 SOR iterations per timestep. 

\subsection{Manufactured Solution} \label{Sec:Incompressible Fluid Flow Problems Manufactured Solution}
We begin with two analytical solutions to Navier-Stokes equations known as the Kovasznay and cylindrical Couette flows. These solutions are useful in error analysis since the exact distributions of the variables are known. The components for velocity and pressure in the Kovasznay flow are given by \cite{nektar_kovasznay}:
\begin{equation}
\begin{split}
u & = 1- \exp(\lambda x) \cos(2 \pi y) \\
v & = \lambda \exp(\lambda x) \sin(2 \pi y) / (2 \pi) \\
p & = p_0 - (\exp(2 \lambda x)/2)
\end{split}
\label{Eq:kovasznay uvp}
\end{equation}
where, $p_0$ is arbitrary reference pressure and $\lambda$ is a parameter defined in terms of Reynolds number:
\begin{equation}
\lambda = \frac{Re}{2} - \left(\frac{Re^2}{4} + 4 \pi^2\right) ^{0.5}
\label{Eq:kovasznay lambda}
\end{equation}
Vertices of the unstructured triangular mesh of the unit square are used as scattered points for the meshless discretization. Three different spatial resolutions are simulated having average $\Delta x$ of $[0.0330,0.0167,0.0101]$.
\par The cylindrical Couette flow consists of an annular region between two cylinders filled with fluid. The inner cylinder is rotating and outer one is kept stationary. This problem has a one dimensional exact solution in radial coordinates. The tangential velocity is given by \cite{white2016fluid}:
\begin{equation}
	v_{\theta}(r) = r_1 \Omega \frac{r_1 r_2}{r_2^2 - r_1^2} \left(\frac{r_2}{r} - \frac{r}{r_2}\right)
	\label{Eq:couette v_theta}
\end{equation}
where, $\Omega$, $r$, $r_1$ and $r_2$ are angular velocity, radial coordinate, inner and outer radii respectively. We solve this problem in Cartesian coordinates and thus, the flow is considered two dimensional. The flow is known to be two dimensional, steady and laminar at low Reynolds number \cite{white2016fluid}. The Reynolds number defined in terms of tangential velocity and inner cylinder's diameter is set to 100. However, the velocity field does not depend on the Reynolds number. We have prescribed the radii and angular velocity to be $r_1=0.5$, $r_2=1$ and $\Omega=2$. As earlier, vertices lying on a triangular mesh in the annular region with an average $\Delta x$ of $[0.0482,0.0357,0.0256]$ are used as the scattered points.

\begin{figure}[H]
	\centering
	\begin{subfigure}[t]{0.45\textwidth}
		\includegraphics[width=\textwidth]{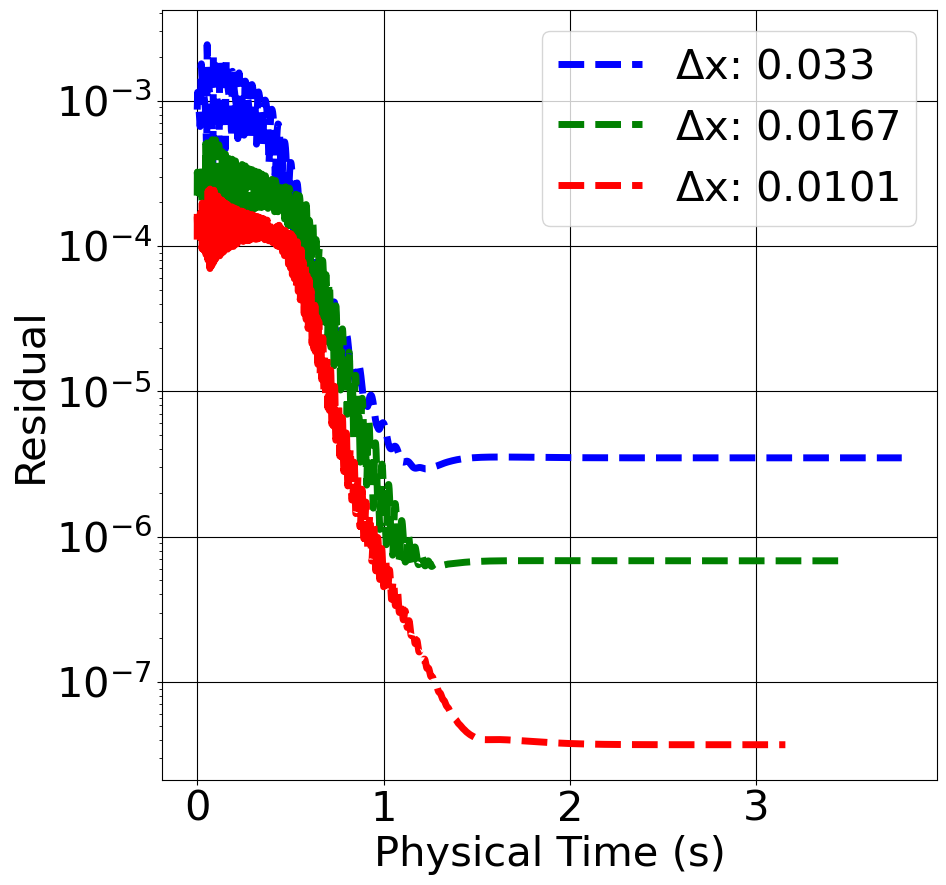}
		\caption{Kovasznay Flow}
	\end{subfigure}
	\hspace{0.05\textwidth}
	\begin{subfigure}[t]{0.45\textwidth}
		\includegraphics[width=\textwidth]{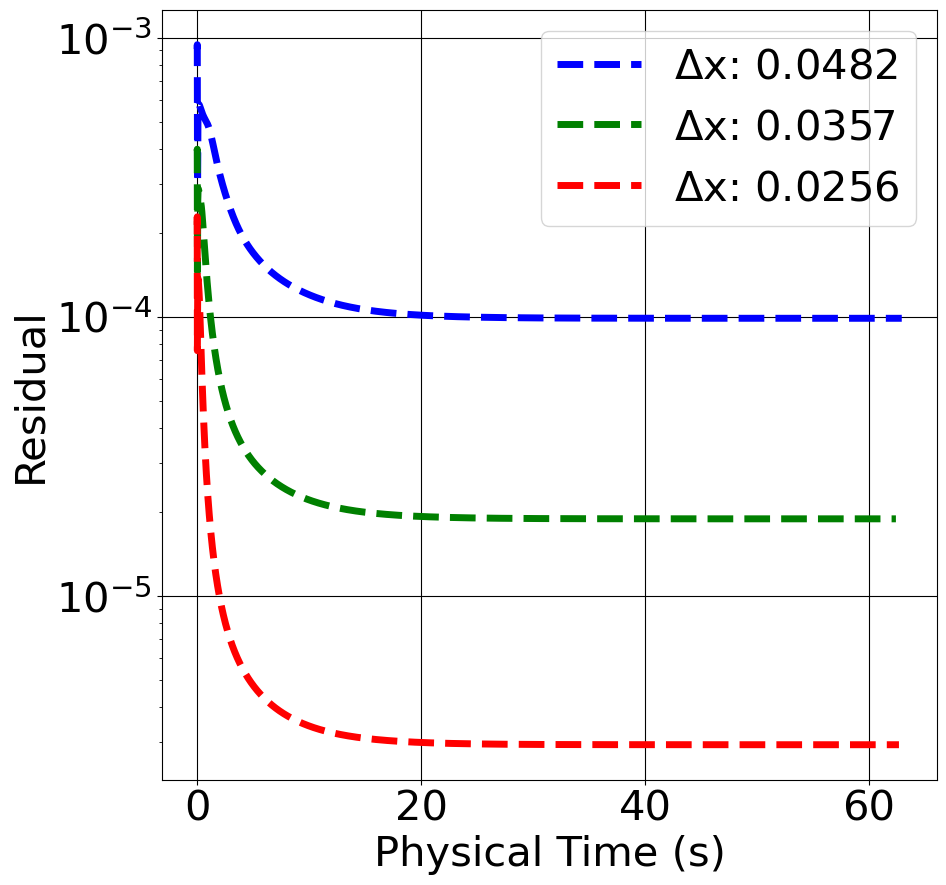}
		\caption{Circular Couette Flow}
	\end{subfigure}
	\caption{Residuals for Varying Grid Sizes for Polynomial Degree of 3}
	\label{Fig:Meshless: kovasznay couette res vs iter polydeg 3}
\end{figure}

\begin{figure}[H]
	\centering
	\begin{subfigure}[t]{0.45\textwidth}
		\includegraphics[width=\textwidth]{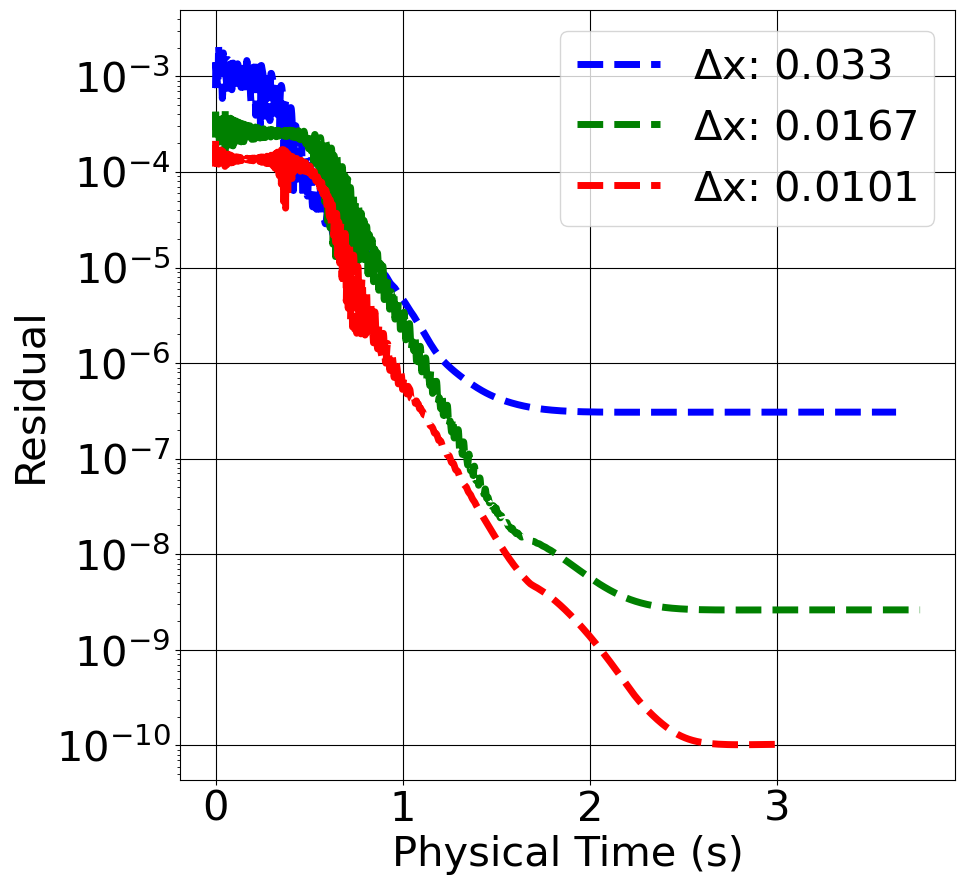}
		\caption{Kovasznay Flow}
	\end{subfigure}
	\hspace{0.05\textwidth}
	\begin{subfigure}[t]{0.45\textwidth}
		\includegraphics[width=\textwidth]{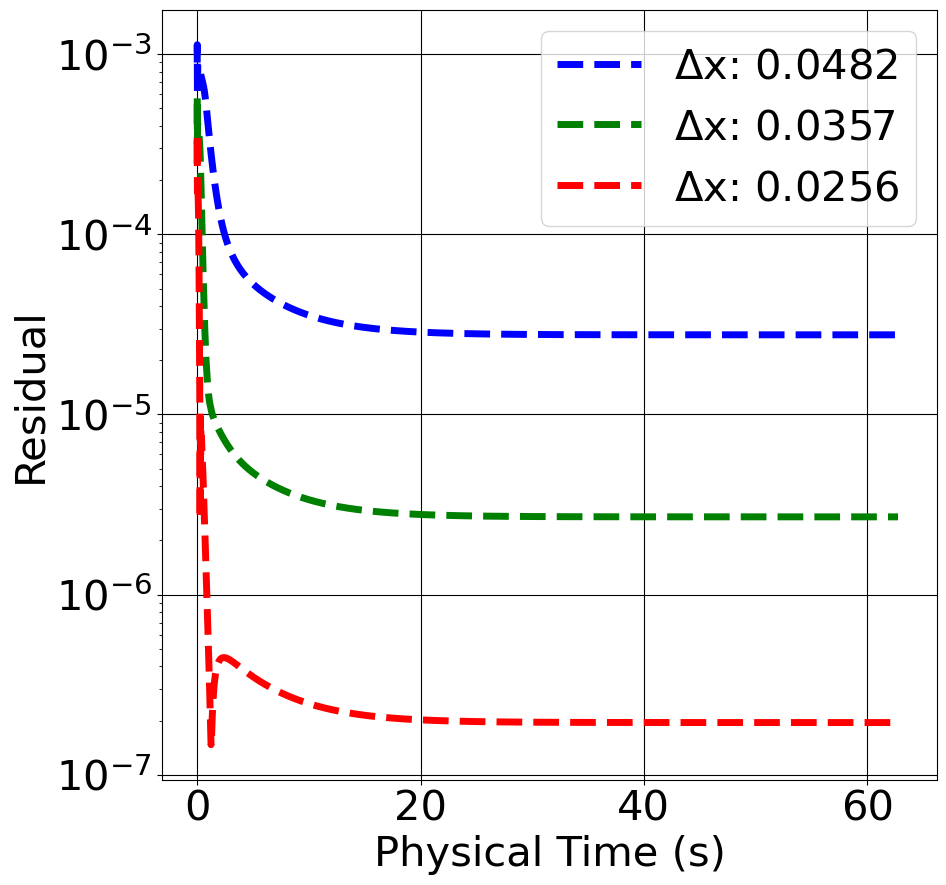}
		\caption{Circular Couette Flow}
	\end{subfigure}
	\caption{Residuals for Varying Grid Sizes for Polynomial Degree of 6}
	\label{Fig:Meshless: kovasznay couette res vs iter polydeg 6}
\end{figure}
\begin{figure}[H]
	\centering
	\begin{subfigure}[t]{0.45\textwidth}
		\includegraphics[width=\textwidth]{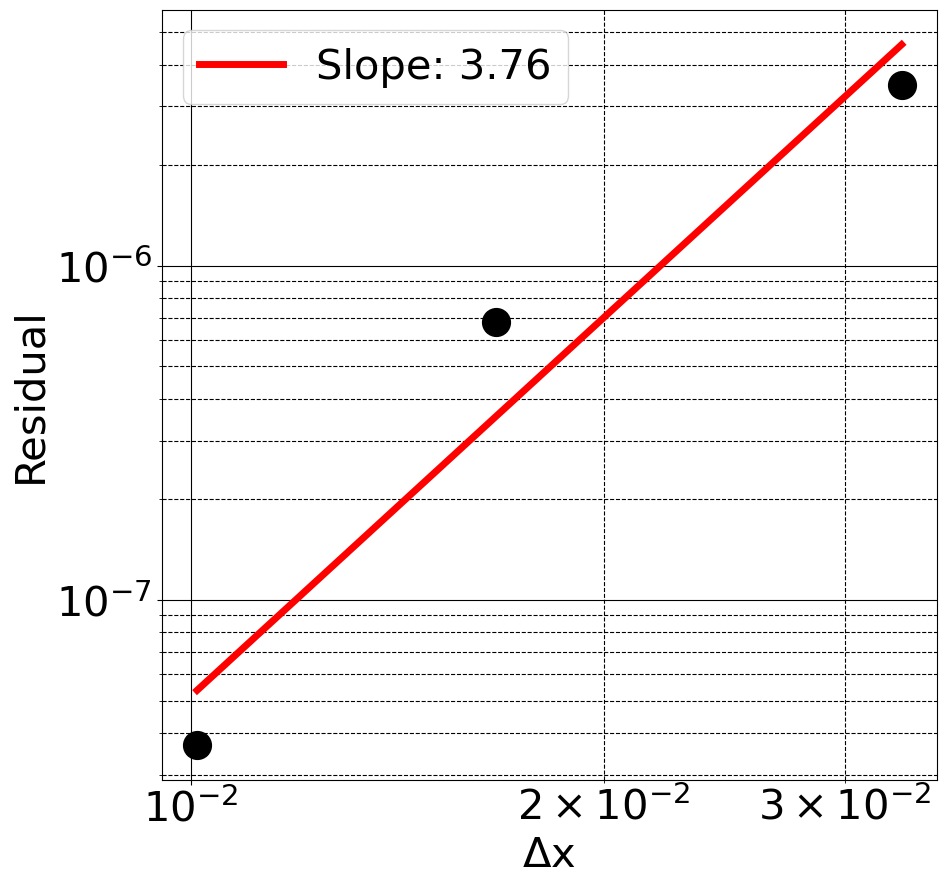}
		\caption{Kovasznay Flow}
	\end{subfigure}
	\hspace{0.05\textwidth}
	\begin{subfigure}[t]{0.45\textwidth}
		\includegraphics[width=\textwidth]{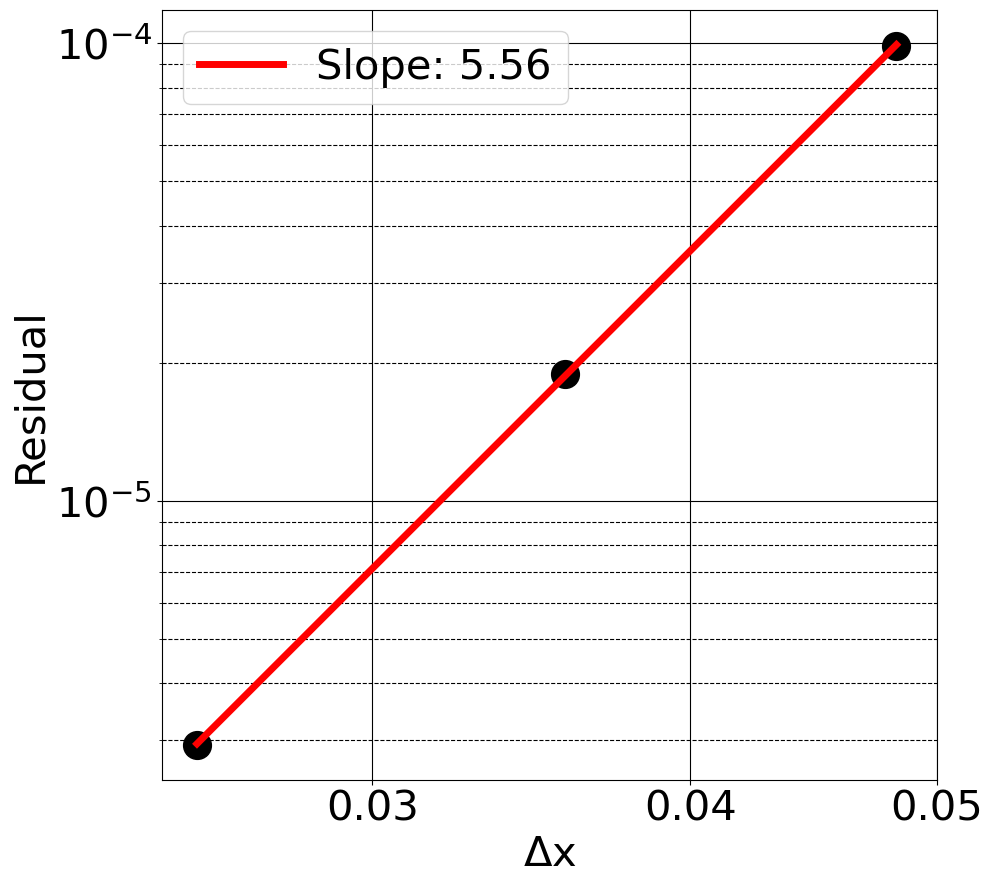}
		\caption{Circular Couette Flow}
	\end{subfigure}
	\caption{Convergence Rate of Final Residuals for Polynomial Degree of 3}
	\label{Fig:Meshless: kovasznay couette rate of conv residual polydeg 3}
\end{figure}
\begin{figure}[H]
	\centering
	\begin{subfigure}[t]{0.45\textwidth}
		\includegraphics[width=\textwidth]{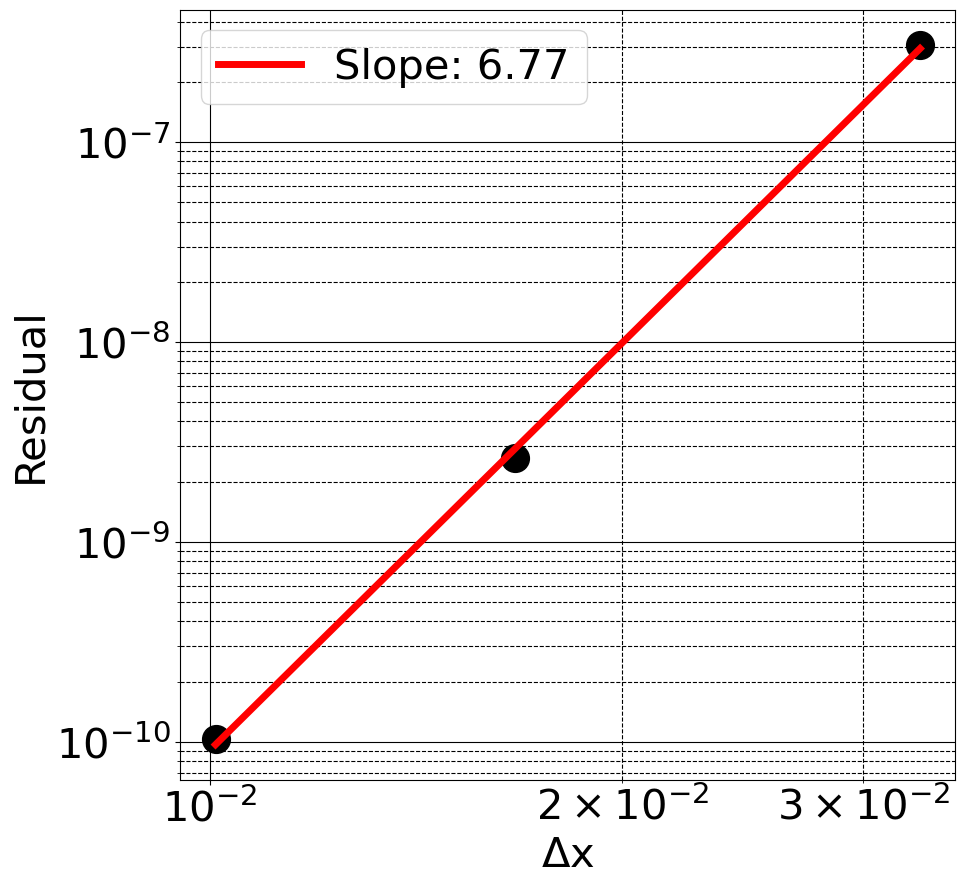}
		\caption{Kovasznay Flow}
	\end{subfigure}
	\hspace{0.05\textwidth}
	\begin{subfigure}[t]{0.45\textwidth}
		\includegraphics[width=\textwidth]{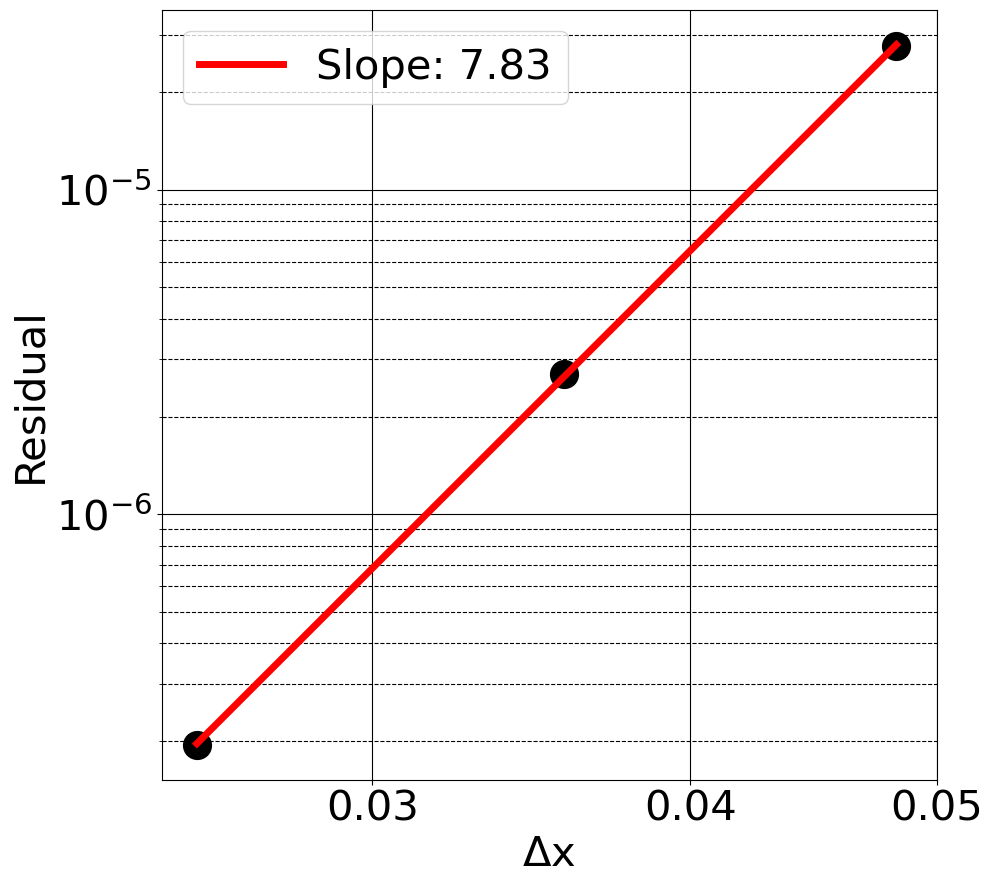}
		\caption{Circular Couette Flow}
	\end{subfigure}
	\caption{Convergence Rate of Final Residuals for Polynomial Degree of 6}
	\label{Fig:Meshless: kovasznay couette rate of conv residual polydeg 6}
\end{figure}

For both the problems, we simulate three different grid resolutions with four degrees of appended polynomials from 3 to 6. As mentioned before, 10-20 SOR iterations are performed at each timestep for the solution of the pressure Poisson equation (\cref{Eq:frac step PPE}). The residual at the end of these iterations is plotted in \cref{Fig:Meshless: kovasznay couette res vs iter polydeg 3,Fig:Meshless: kovasznay couette res vs iter polydeg 6} as a function of physical time. We have only plotted the cases of polynomial degrees 3 and 6. Other polynomial degrees show similar trends but higher values of final residuals than the polynomial degree of 6. Similar to the examples in \cref{Sec:Poisson Equation with Manufactured Solution Meshless Discretization}, the final stationary value reduces with the truncation error for larger number of scattered points or higher polynomial degree. \Cref{Fig:Meshless: kovasznay couette rate of conv residual polydeg 3,Fig:Meshless: kovasznay couette rate of conv residual polydeg 6} plot the value of the final residual with average $\Delta x$ on a logarithmic scale. The best fit lines show that the rates of convergence are higher than the respective polynomial degrees.

\begin{figure}[H]
	\centering
	\includegraphics[width=0.6\textwidth]{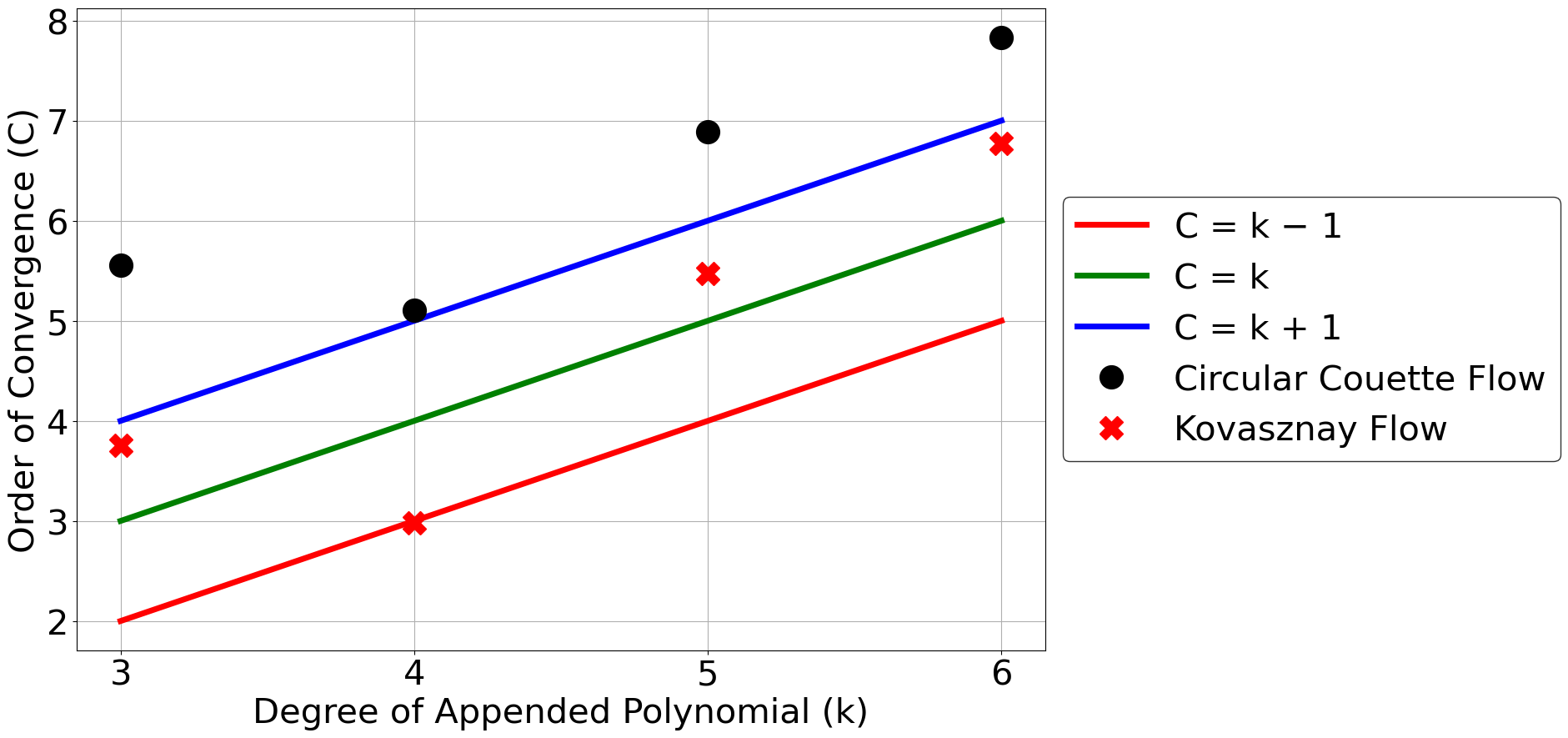}
	\caption{Composite Convergence Rate of Final Residuals with Degrees of Appended Polynomial}
	\label{Fig:Meshless: kovasznay couette rate of conv residual with polydeg}
\end{figure}

\begin{figure}[H]
	\centering
	\includegraphics[width=0.6\textwidth]{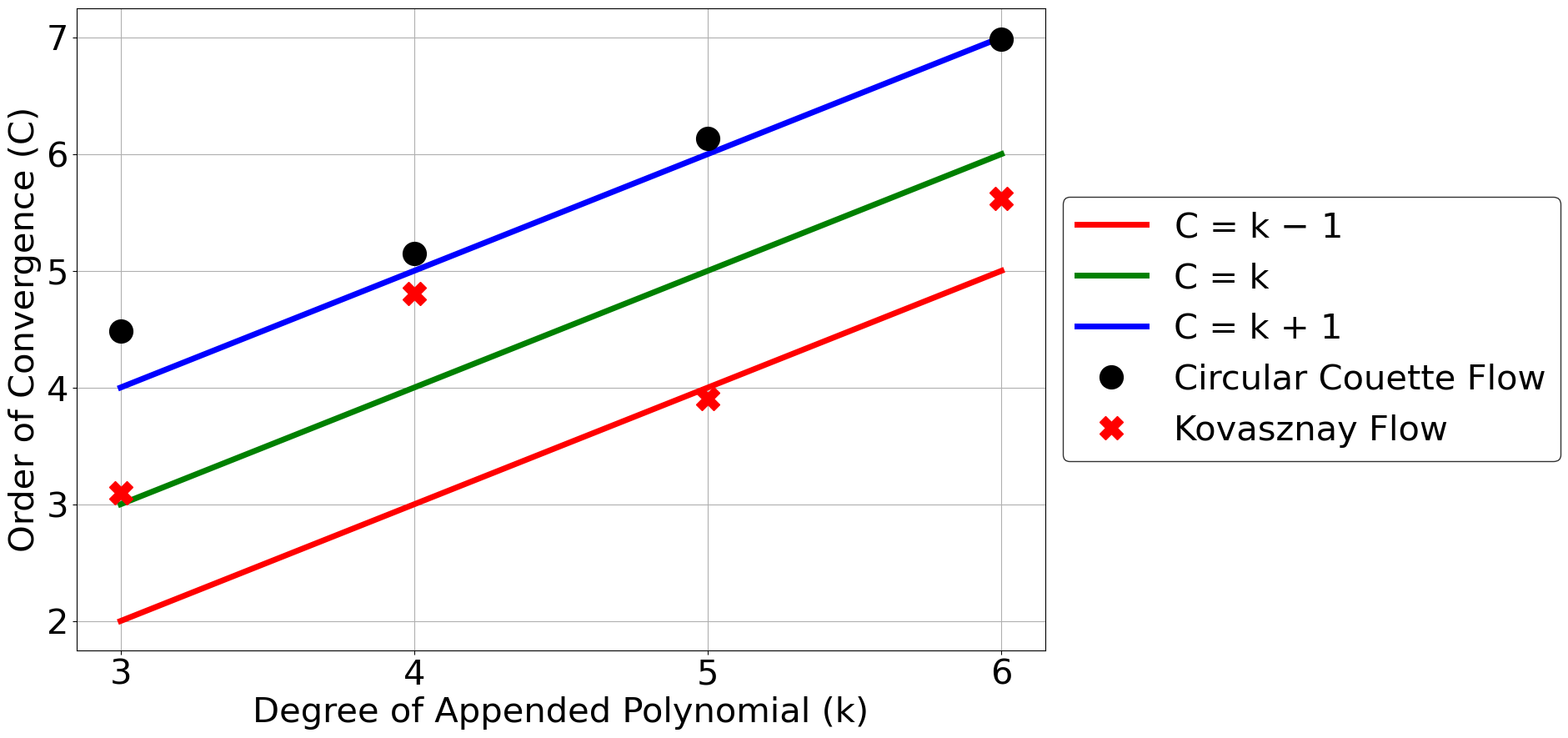}
	\caption{Composite Convergence Rate of Solution Error with Degrees of Appended Polynomial}
	\label{Fig:Meshless: kovasznay couette rate of conv error with polydeg}
\end{figure}

Similar analyses are performed for other polynomial degrees and the rates of convergence for the final residual are plotted in a composite plot in \cref{Fig:Meshless: kovasznay couette rate of conv residual with polydeg} with polynomial degree. These rates are plotted for both the Kovasznay and Couette flow with separate markers. From the reference lines, we can see that the rate of convergence increases with polynomial degree and is typically higher than the degree of the appended polynomial.
We now present the L$_1$ norm of error in the converged velocities and its rate of convergence with increasing number of points and polynomial degrees. Since the exact solution is known, the error in both the velocity components is defined in a manner similar to \cref{Eq:rel res error norms}. Because the Navier-Stokes equations have Laplacian and gradient operators, it is expected that the PHS-RBF stencils with appended polynomials of degree $k$ converge at rates between $\mathcal{O}(k)$ and $\mathcal{O}(k-1)$ \cite{shahane2021high,flyer2016onrole_I}. We find that is indeed the case, as most of the rates are found to follow these reference lines. However, \cref{Fig:Meshless: kovasznay couette rate of conv residual with polydeg,Fig:Meshless: kovasznay couette rate of conv error with polydeg} have a few outliers probably because the value of $\Delta x$ used for the best fit line is only an approximate measure of spatial resolution.

\subsection{Analysis of Complex Flows} \label{Sec:Incompressible Fluid Flow Problems Analysis of Complex Flows}
In this section, we further consider two complex fluid flow problems which have applications in heat exchangers and manufacturing. Both the problems have complex recirculation zones. Similar to previous examples, we use time marching with the fractional step method \cite{shahane2021high} to reach steady state. The pressure Poisson equation is solved by SOR iterations. As before, convergence of the stationary value of the residual is studied.
\subsubsection{Steady Flow in a Periodic Converging-Diverging Passage}
\begin{figure}[H]
	\centering
	\includegraphics[width=0.6\textwidth]{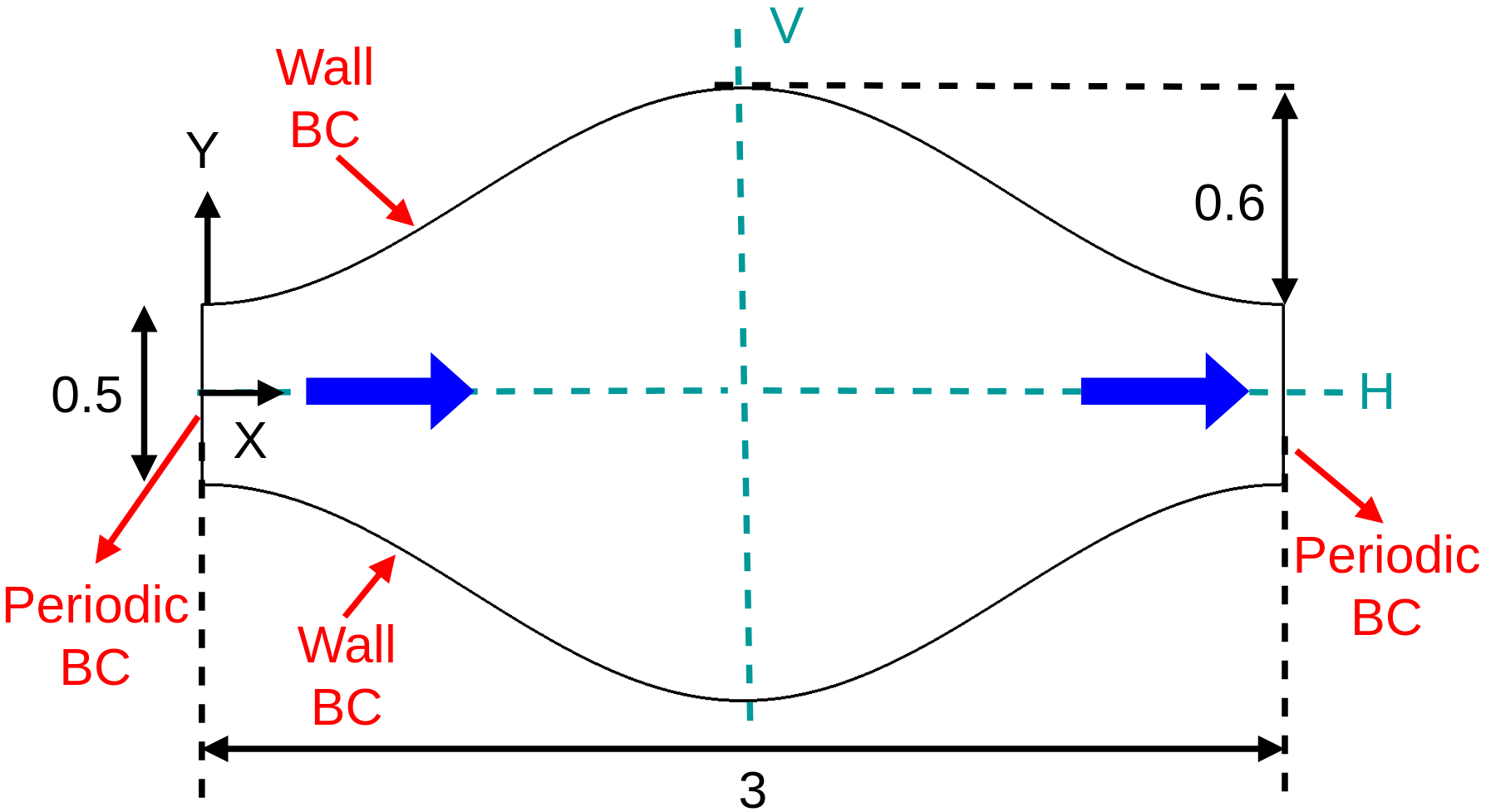}
	\caption{Geometry and Boundary Conditions of the Periodic Converging-Diverging Passage}
	\label{Fig:Meshless: periodic bellow schematic}
\end{figure}
First we consider a converging-diverging channel showed in \cref{Fig:Meshless: periodic bellow schematic}. The left and right boundaries of the domain are set to periodic condition. The bottom and top boundaries are prescribed as walls with no-slip and no-penetration. For details of the implementation of periodicity in meshless method, please refer to our previous work \cite{shahane2021semi}. We define the minimum width of the channel ($L=0.5$) as the characteristic length and the density ($\rho$) is set to unity. A mean pressure gradient is prescribed in X direction to drive the fluid flow:
\begin{equation}
\frac{\partial \bar{p}}{\partial x} = \frac{12 \mu}{L^2}
\label{Eq:periodic bellow dp_dx}
\end{equation}
After the steady state is obtained, the velocity component in X direction is averaged over the vertical direction at the left boundary to estimate the flow Reynolds number: $Re=\rho \bar{u} L / \mu$. The velocity component and pressure are initialized to be zero. The viscosity is set to 0.01 which yields a steady state flow with a Reynolds number of 111.3. The velocity components and vorticity are interpolated to a uniform array of points along the marked horizontal (H) and vertical (V) lines (shown in \cref{Fig:Meshless: periodic bellow schematic}) to study grid independency of the results.
\begin{figure}[H]
	\centering
	\begin{subfigure}[t]{0.45\textwidth}
		\includegraphics[width=\textwidth]{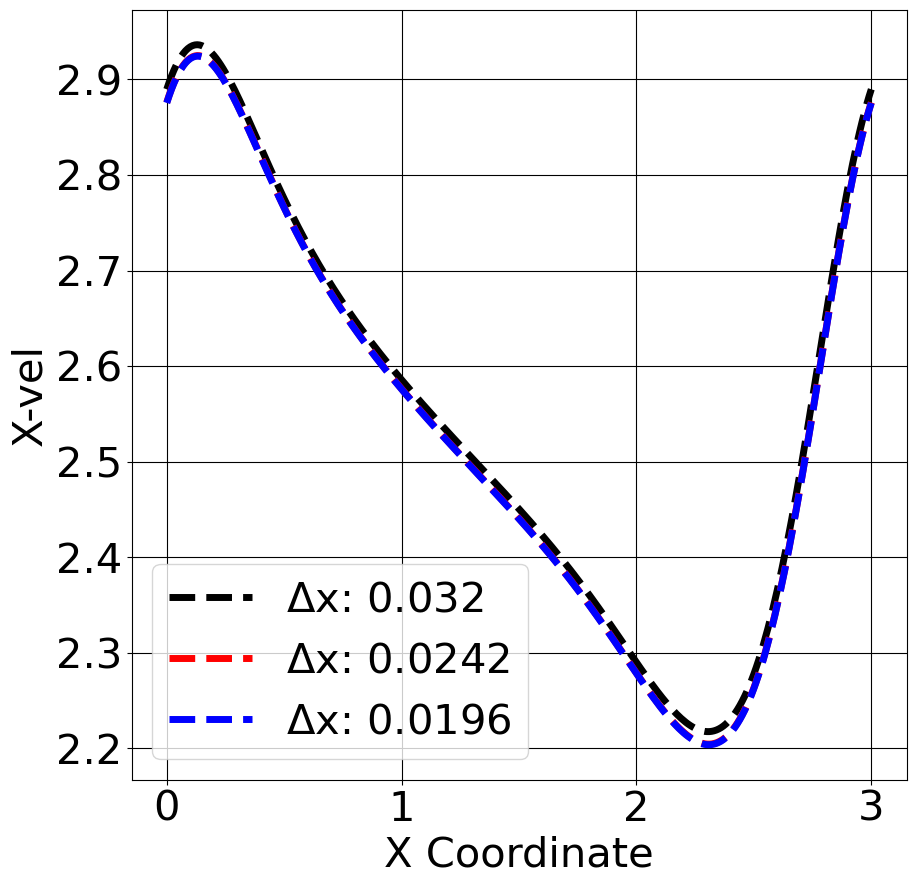}
		\caption{X-Vel: Horizontal Center Line (H)}
		\label{Fig:Meshless: periodic bellow Grid Independence X-vel H}
	\end{subfigure}
	\hspace{0.05\textwidth}
	\begin{subfigure}[t]{0.45\textwidth}
		\includegraphics[width=\textwidth]{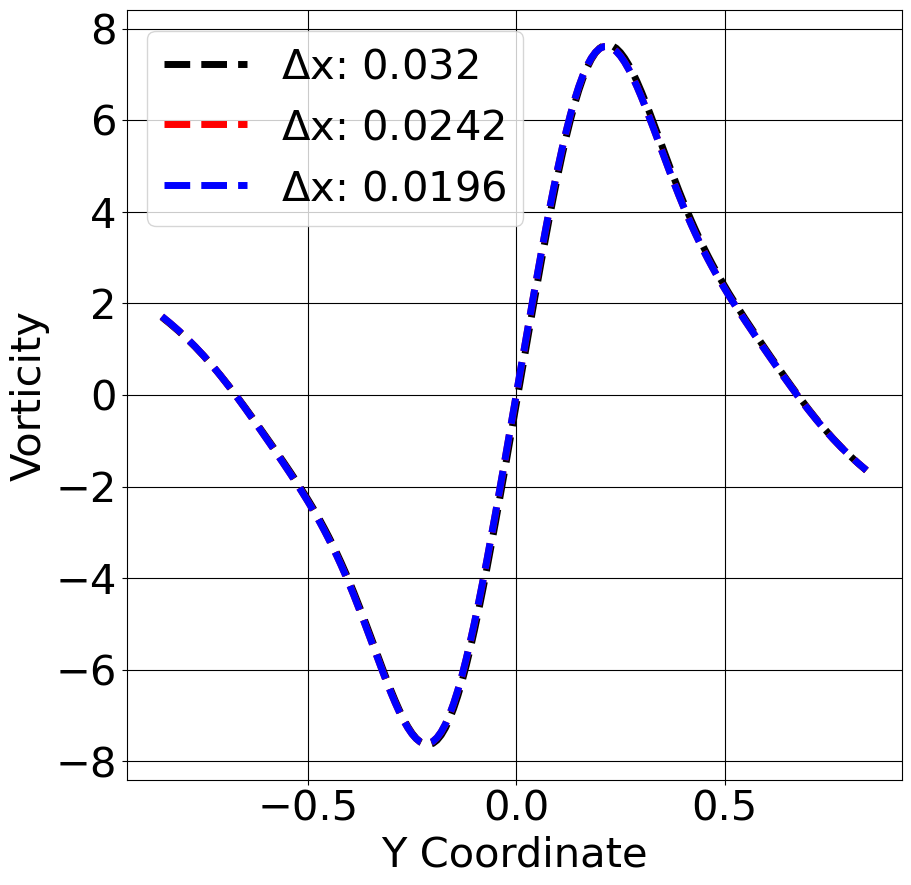}
		\caption{Vorticity: Vertical Center Line (V)}
		\label{Fig:Meshless: periodic bellow Grid Independence vorticity V}
	\end{subfigure}
	\caption{Grid Independence for Re $=111.3$ with Polynomial Degree 6}
	\label{Fig:Meshless: periodic bellow Grid Independence}
\end{figure}

\begin{figure}[H]
	\centering
	\begin{subfigure}[t]{0.45\textwidth}
		\includegraphics[width=\textwidth]{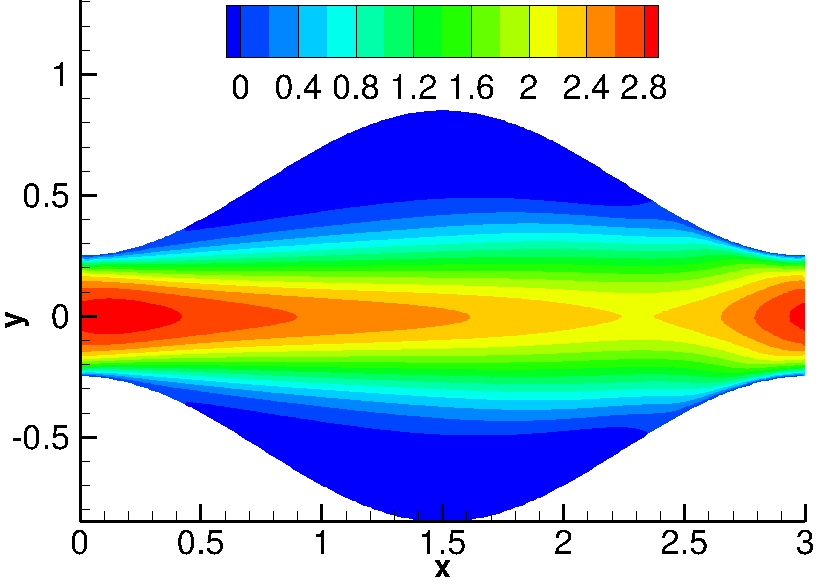}
		\caption{X-Vel}
	\end{subfigure}
	\hspace{0.05\textwidth}
	\begin{subfigure}[t]{0.45\textwidth}
		\includegraphics[width=\textwidth]{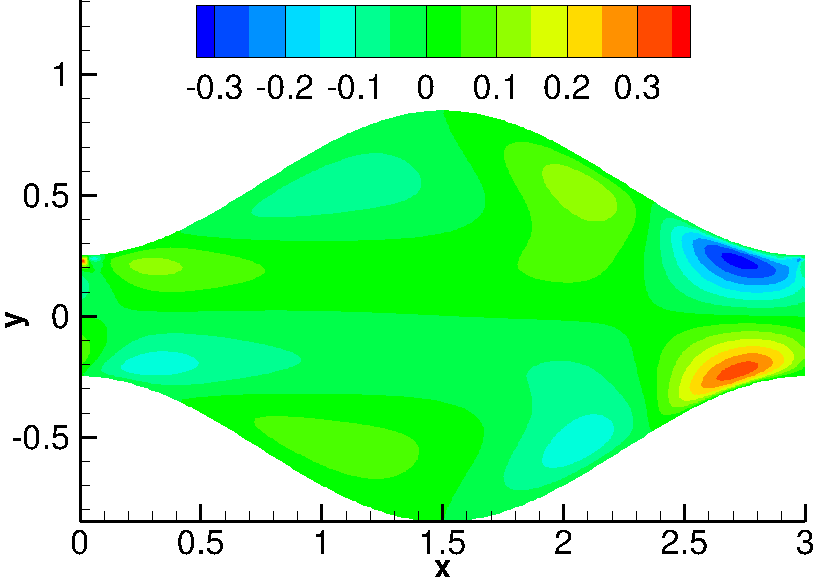}
		\caption{Y-Vel}
	\end{subfigure}\vspace{0.25cm}
	\begin{subfigure}[t]{0.45\textwidth}
		\includegraphics[width=\textwidth]{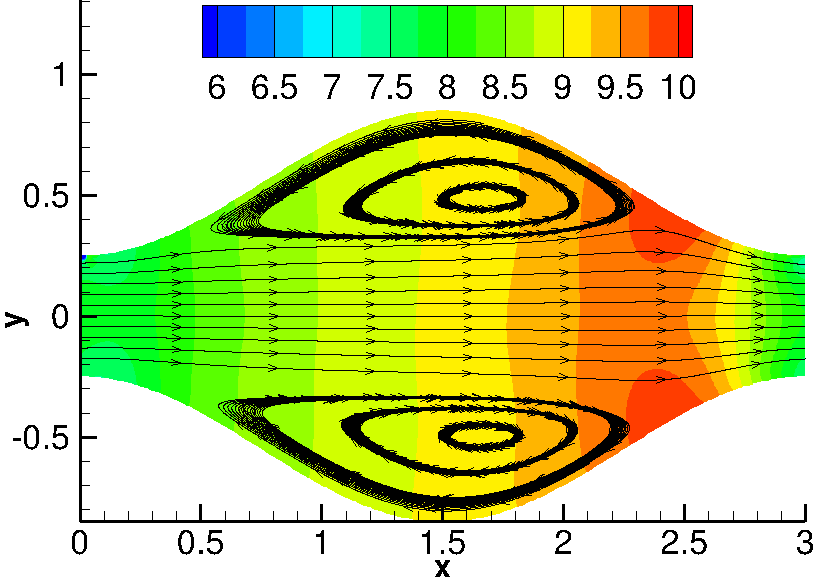}
		\caption{Pressure Variations and Streamlines}
		\label{Fig:Meshless: periodic bellow contours p str}
	\end{subfigure}
	\caption{Contour Plots for Re $=111.3$ ($\Delta x=0.0196$ with Polynomial Degree 6)}
	\label{Fig:Meshless: periodic bellow contours}
\end{figure}

Three different unstructured point distributions are used with average $\Delta x$ of $[0.0320,0.0242,0.0196]$ and the polynomial degree is varied from 3 to 6. \Cref{Fig:Meshless: periodic bellow Grid Independence} plots the velocity and vorticity interpolated to H and V lines for the three point distributions and polynomial degree of 6. In \cref{Fig:Meshless: periodic bellow Grid Independence X-vel H}, there is a slight mismatch for the coarsest grid but the results of the two finer grids overlap. On the other hand, \cref{Fig:Meshless: periodic bellow Grid Independence vorticity V} shows that the vorticity estimates overlap for all the three grids. Thus, grid independent solutions are obtained. Contours estimated using the finest grid and polynomial degree of 6 are plotted in \cref{Fig:Meshless: periodic bellow contours}. We obtain a steady state solution for the Reynolds number of 111.3. The contours plotted in \cref{Fig:Meshless: periodic bellow contours p str} are pressure variations over the mean pressure imposed to drive the periodic fully-developed flow. As expected, we observe that the steady flow is symmetric about the horizontal center line.

\begin{figure}[H]
	\centering
	\begin{subfigure}[t]{0.45\textwidth}
		\includegraphics[width=\textwidth]{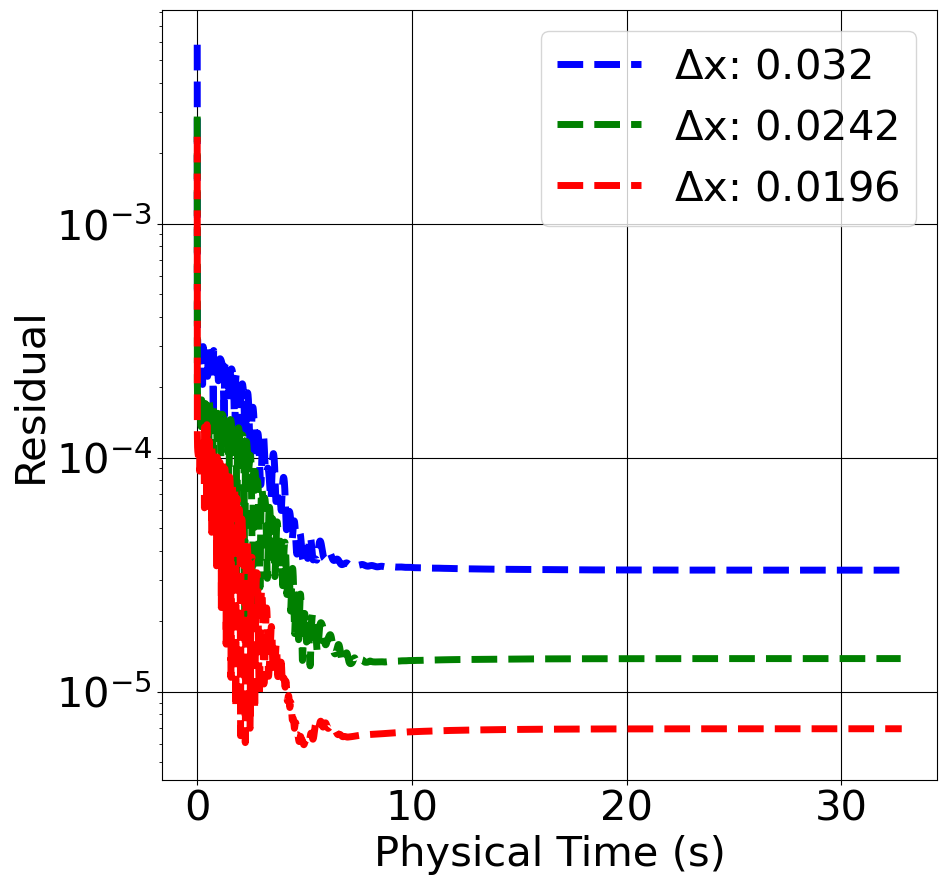}
		\caption{Degree of Appended Polynomial: 3}
		\label{Fig:Meshless: periodic bellow res vs iter polydeg 3}
	\end{subfigure}
	\hspace{0.05\textwidth}
	\begin{subfigure}[t]{0.45\textwidth}
		\includegraphics[width=\textwidth]{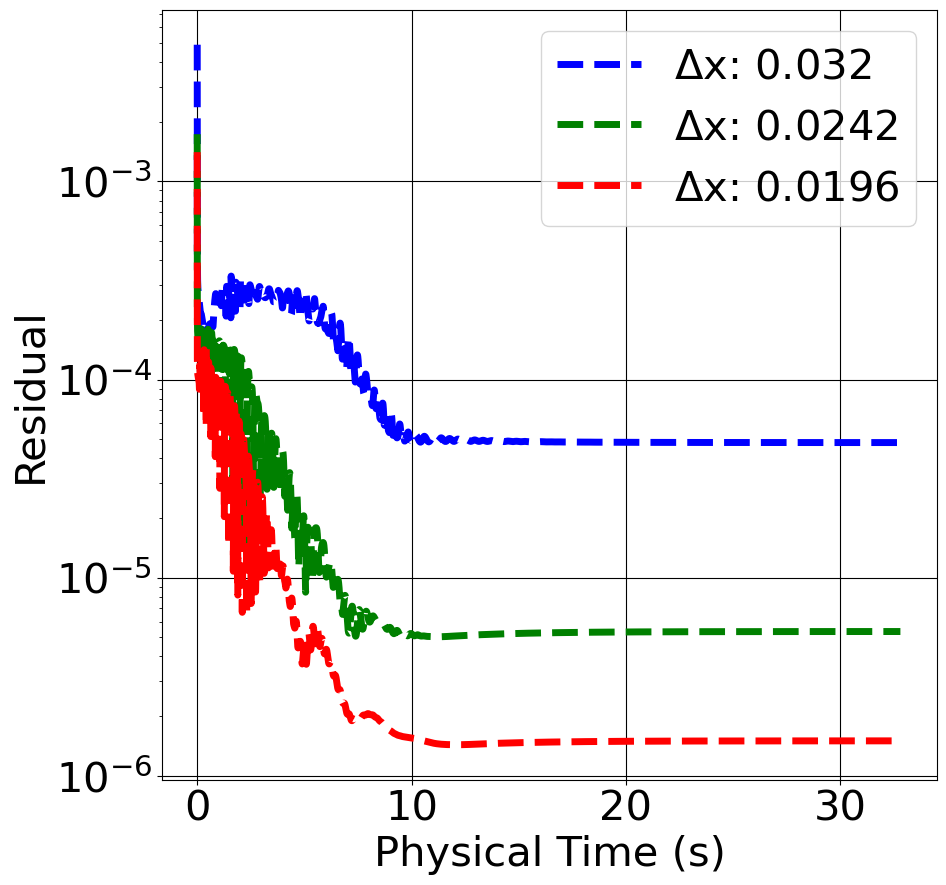}
		\caption{Degree of Appended Polynomial: 6}
		\label{Fig:Meshless: periodic bellow res vs iter polydeg 6}
	\end{subfigure}
	\caption{Residuals for Varying Grid Sizes}
	\label{Fig:Meshless: periodic bellow res vs iter}
\end{figure}

\begin{figure}[H]
	\centering
	\begin{subfigure}[t]{0.45\textwidth}
		\includegraphics[width=\textwidth]{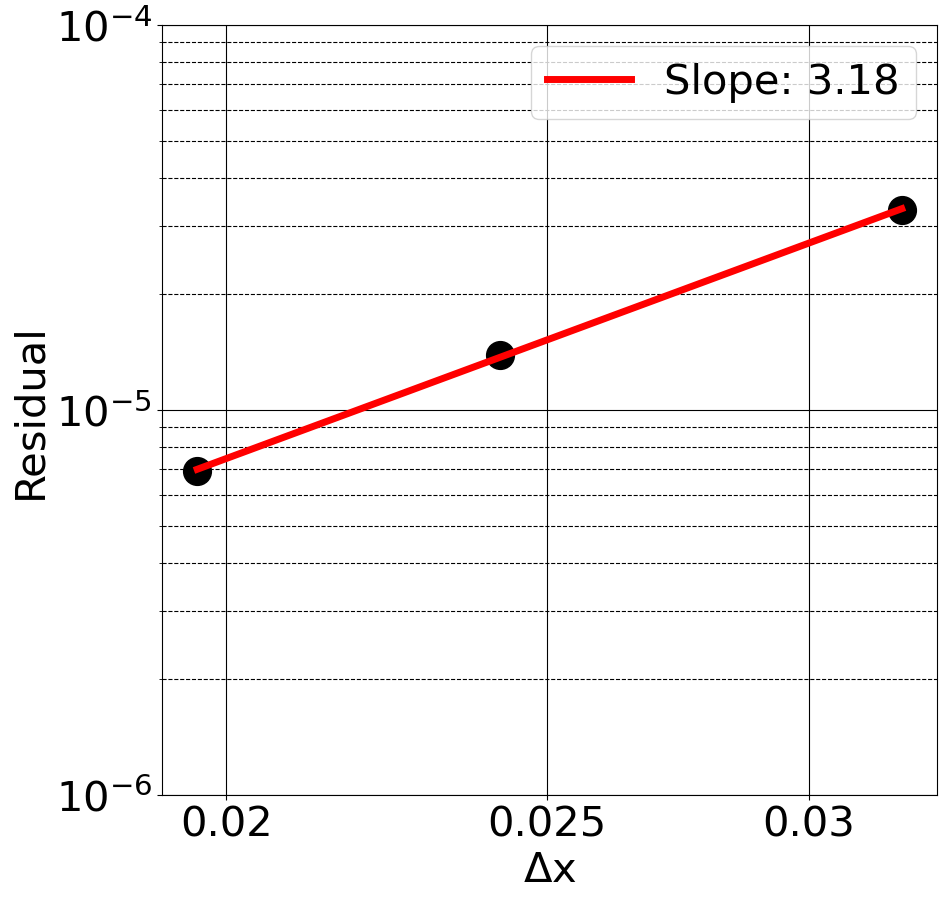}
		\caption{Degree of Appended Polynomial: 3}
	\end{subfigure}
	\hspace{0.05\textwidth}
	\begin{subfigure}[t]{0.45\textwidth}
		\includegraphics[width=\textwidth]{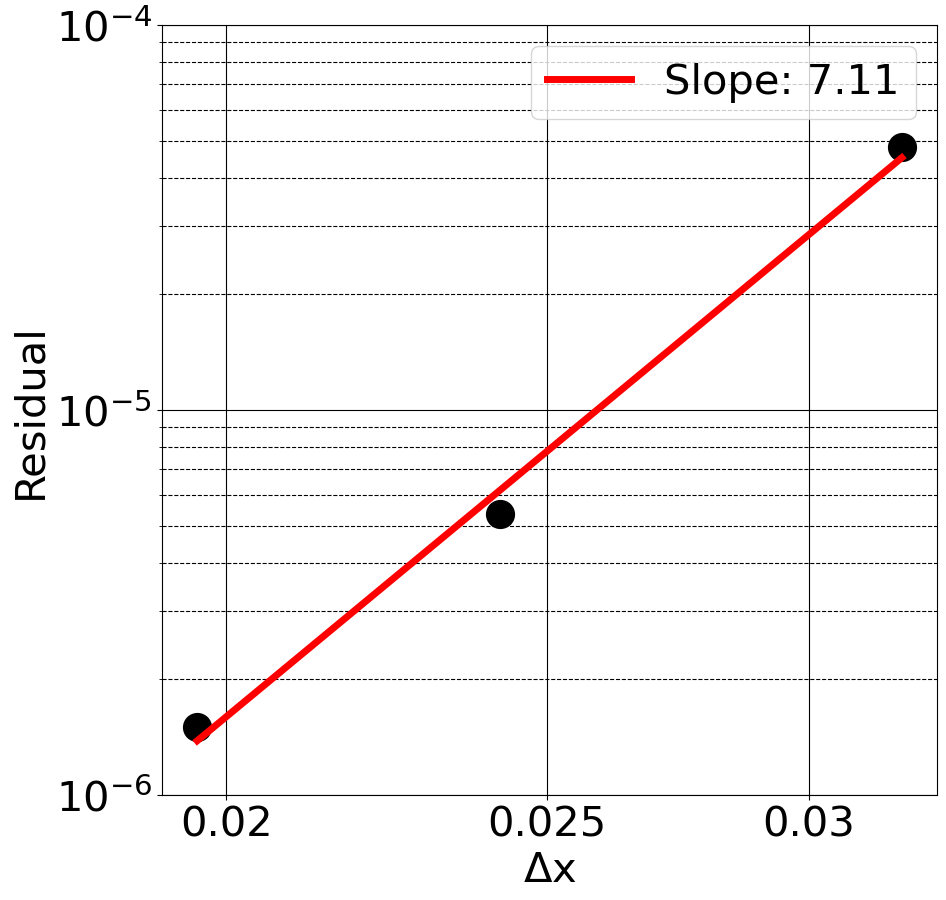}
		\caption{Degree of Appended Polynomial: 6}
	\end{subfigure}
	\caption{Convergence Rate of Final Residuals}
	\label{Fig:Meshless: periodic bellow rate of conv residual}
\end{figure}

Now we analyze the residual of pressure Poisson equation (\cref{Eq:frac step PPE}) as a function of physical time. \Cref{Fig:Meshless: periodic bellow res vs iter} plots the residual at the end of each timestep for three grid resolutions and polynomial degrees of 3 and 6. Similar to previous cases, the residual stabilizes after initial timesteps. We observe that refining the grid or increasing the polynomial degree helps in reducing the stationary value of the residual which in turn increases the solution accuracy. Moreover, \cref{Fig:Meshless: periodic bellow res vs iter polydeg 6} shows that the gap between the three lines is wider than the case of \cref{Fig:Meshless: periodic bellow res vs iter polydeg 3}. This is an indication of higher order of convergence obtained due to higher degree of the appended polynomial. This effect can be further observed in \cref{Fig:Meshless: periodic bellow rate of conv residual} which plots the stationary value of the residual with average $\Delta x$ on a logarithmic scale. The best fit line estimates the slopes which indicate the rates of convergence to be 3.18 and 7.11 for the polynomial degrees of the 3 and 6 respectively. A plot summarizing rates of convergence is plotted and discussed later (\cref{Fig:Meshless: applications rate of conv residual with polydeg}).

\subsubsection{Flow Induced by Two Rotating Cylinders}
\begin{figure}[H]
	\centering
	\includegraphics[width=0.6\textwidth]{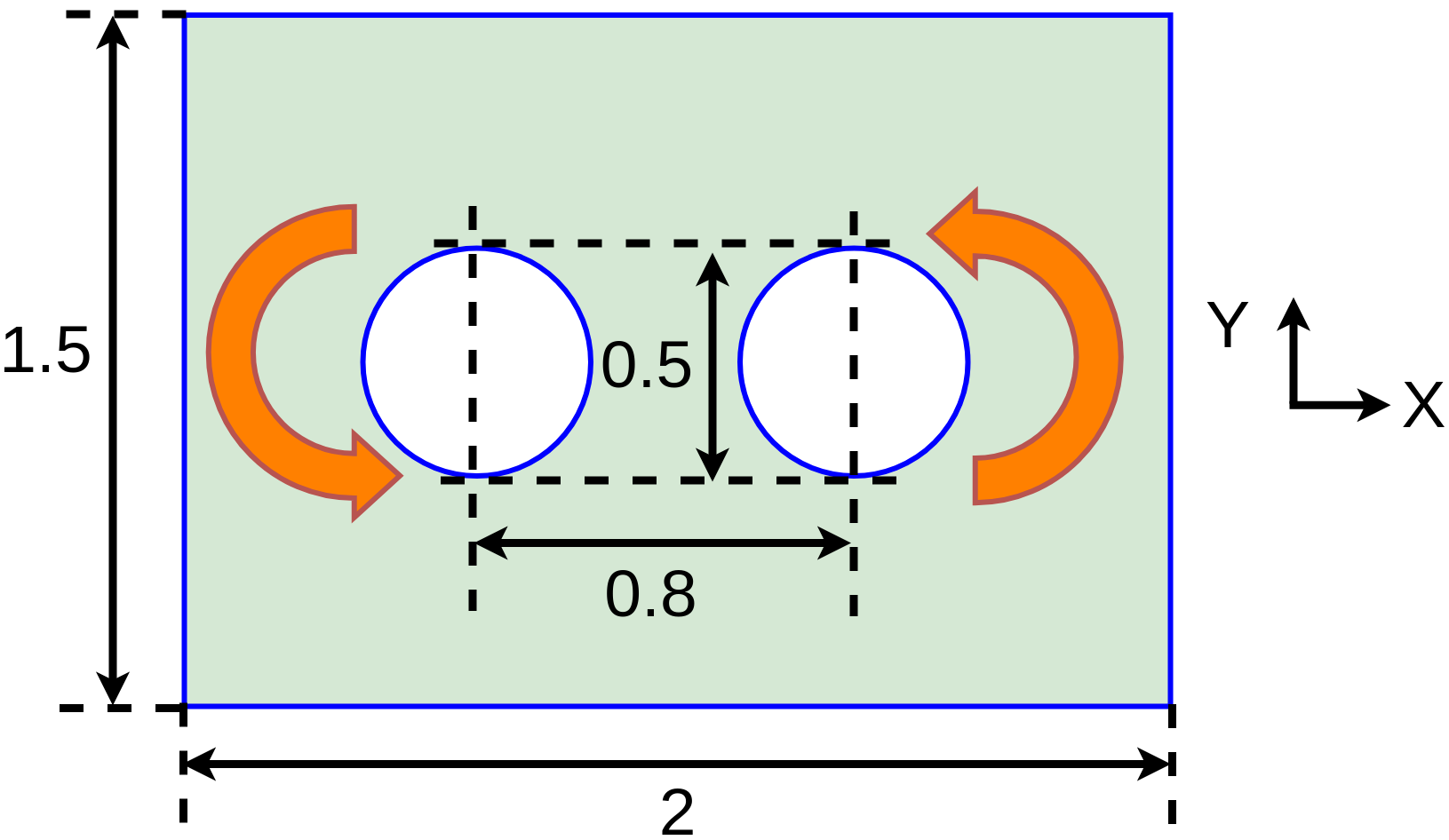}
	\caption{Schematic}
	\label{Fig:Meshless: double circ rect schematic}
\end{figure}
We now consider a complex fluid flow generated by two rotating circular cylinders inside a stationary rectangular enclosure. Region between the cylinders and the rectangle is filled with initially stationary fluid. Both the cylinders rotate in counter-clockwise direction with a angular velocity of 1 radian per second. \Cref{Fig:Meshless: double circ rect schematic} shows the dimensions of the domain. The Reynolds number is based on the cylinder radius and its rotational velocity: $Re=\rho \omega r_i^2 /\mu$ where, $\omega=1$, $\rho=1$, $r_i=0.25$ and $\mu$ is computed from the Reynolds number. The Reynolds number is set to 100, for which we observe a steady two dimensional flow.

\begin{figure}[H]
	\centering
	\begin{subfigure}[t]{0.45\textwidth}
		\includegraphics[width=\textwidth]{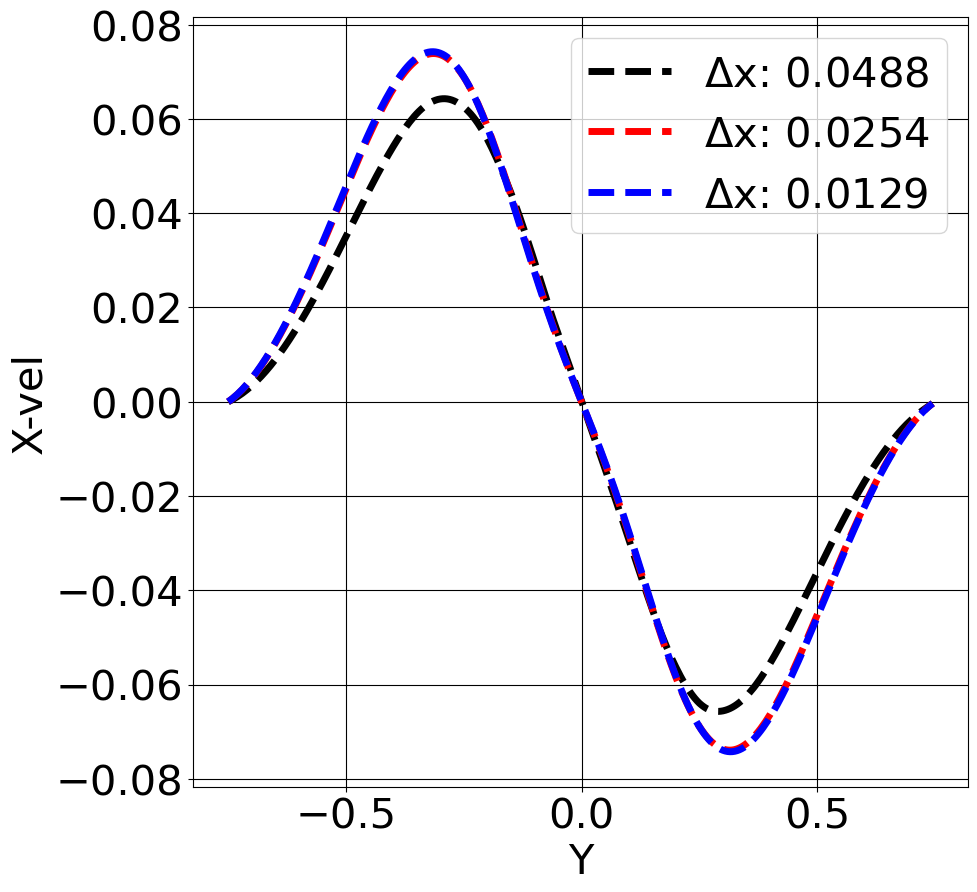}
		\caption{X-Vel on Vertical Center Line}
	\end{subfigure}
	\hspace{0.05\textwidth}
	\begin{subfigure}[t]{0.45\textwidth}
		\includegraphics[width=\textwidth]{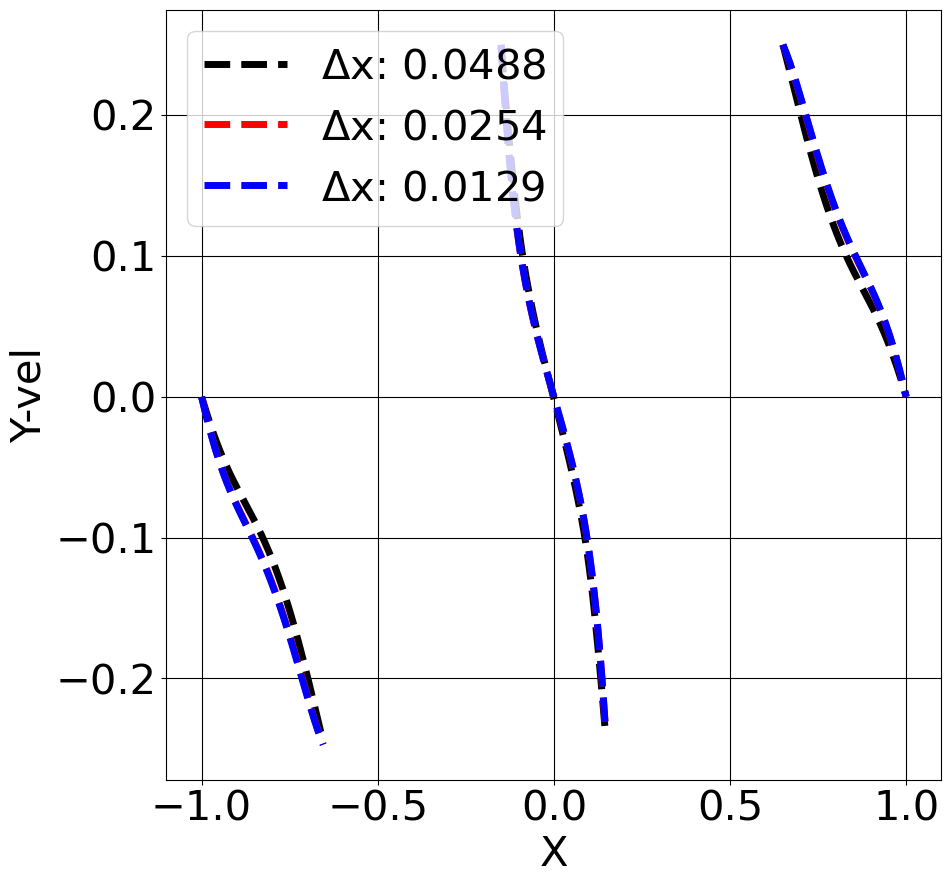}
		\caption{Y-Vel on Horizontal Center Line}
	\end{subfigure}
	\caption{Grid Independence for Re $=100$ with Polynomial Degree 6}
	\label{Fig:Meshless: double circ rect Grid Independence}
\end{figure}
We have simulated this flow with three different point distributions corresponding to average $\Delta x$ of $[0.0488,0.0254,0.0129]$. In order to analyze the effect of point spacing on the flow field, the velocity components are interpolated along horizontal and vertical center lines using the PHS-RBF interpolation. \Cref{Fig:Meshless: double circ rect Grid Independence} show grid independence results for a polynomial degree of 6. It can be seen that though there is a slight mismatch of the coarsest grid, the results from the two finer grids overlap.
\begin{figure}[H]
	\centering
	\begin{subfigure}[t]{0.45\textwidth}
		\includegraphics[width=\textwidth]{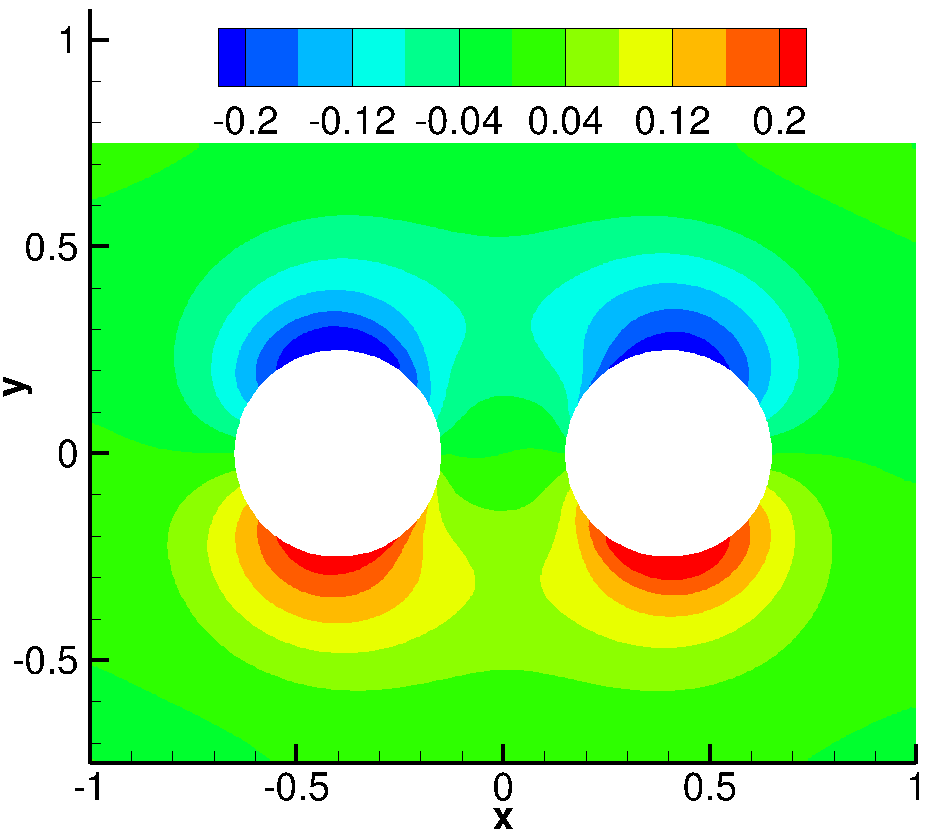}
		\caption{X-Vel}
	\end{subfigure}
	\hspace{0.05\textwidth}
	\begin{subfigure}[t]{0.45\textwidth}
		\includegraphics[width=\textwidth]{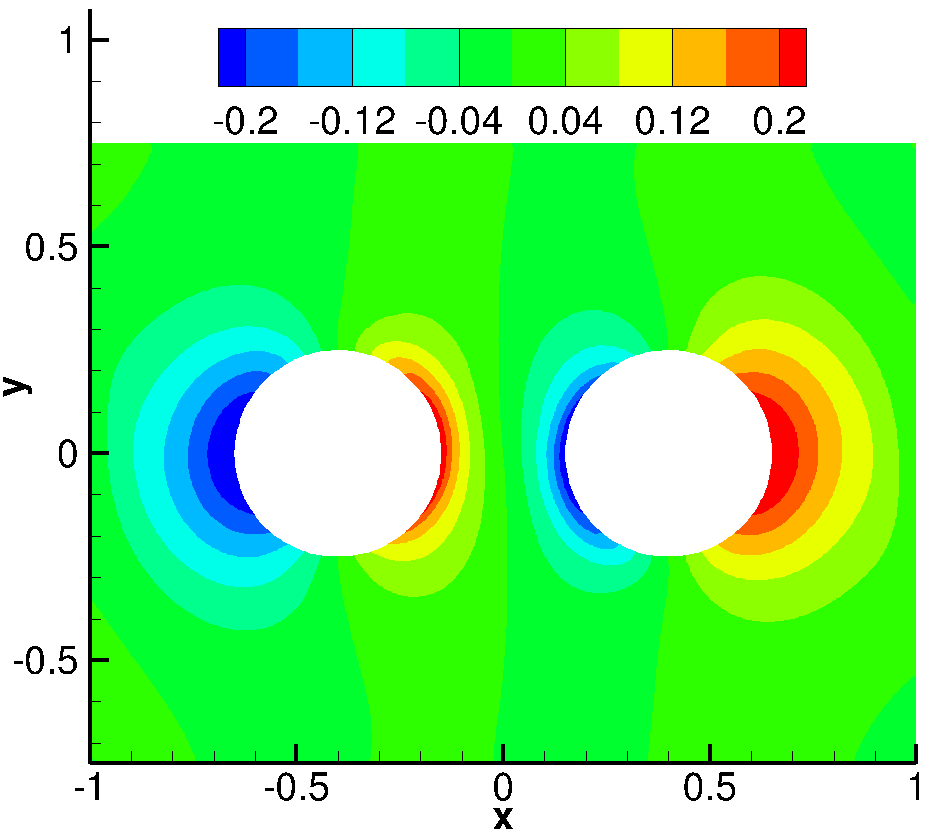}
		\caption{Y-Vel}
	\end{subfigure}\vspace{0.25cm}
	\begin{subfigure}[t]{0.45\textwidth}
		\includegraphics[width=\textwidth]{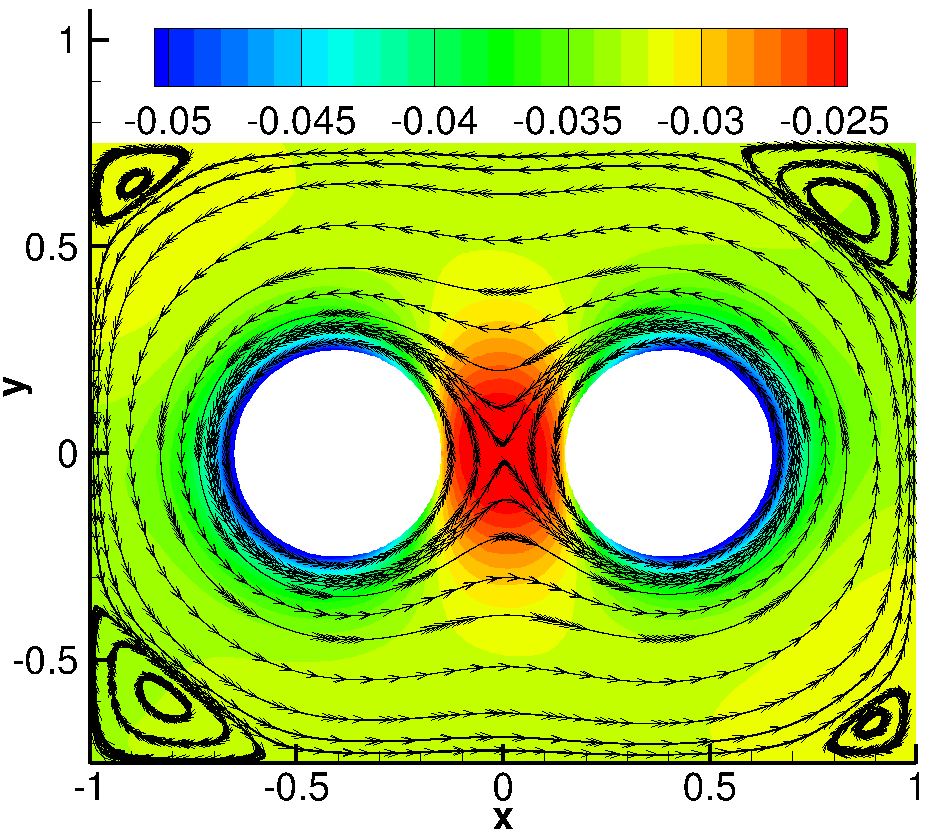}
		\caption{Pressure Variations and Streamlines}
	\end{subfigure}
	\caption{Contour Plots for Re $=100$ ($\Delta x=0.0129$ with Degree of Appended Polynomial: 6)}
	\label{Fig:Meshless: double circ rect contours}
\end{figure}
\Cref{Fig:Meshless: double circ rect contours} plots the contours of velocity components and pressure superposed with streamlines. Due to the destructive interference of the velocities, there is a formation of a stagnation point in the region between the two cylinders. Thus, the pressure at that point is highest and velocity is zero. Moreover, two pairs of secondary vortices are also formed at the corners of the rectangle. The vortices in the first and third quadrant are larger than those in the second and fourth quadrant due to the counter-clockwise rotations of the cylinders.

\begin{figure}[H]
	\centering
	\begin{subfigure}[t]{0.45\textwidth}
		\includegraphics[width=\textwidth]{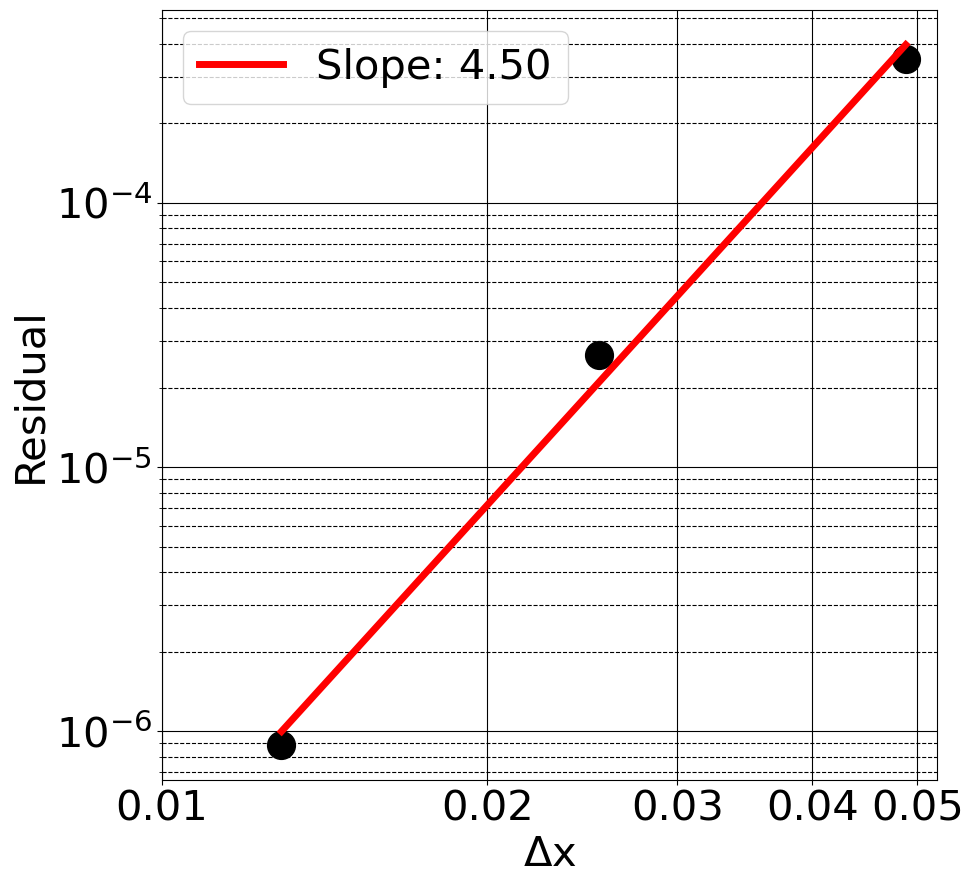}
		\caption{Degree of Appended Polynomial: 3}
	\end{subfigure}
	\hspace{0.05\textwidth}
	\begin{subfigure}[t]{0.45\textwidth}
		\includegraphics[width=\textwidth]{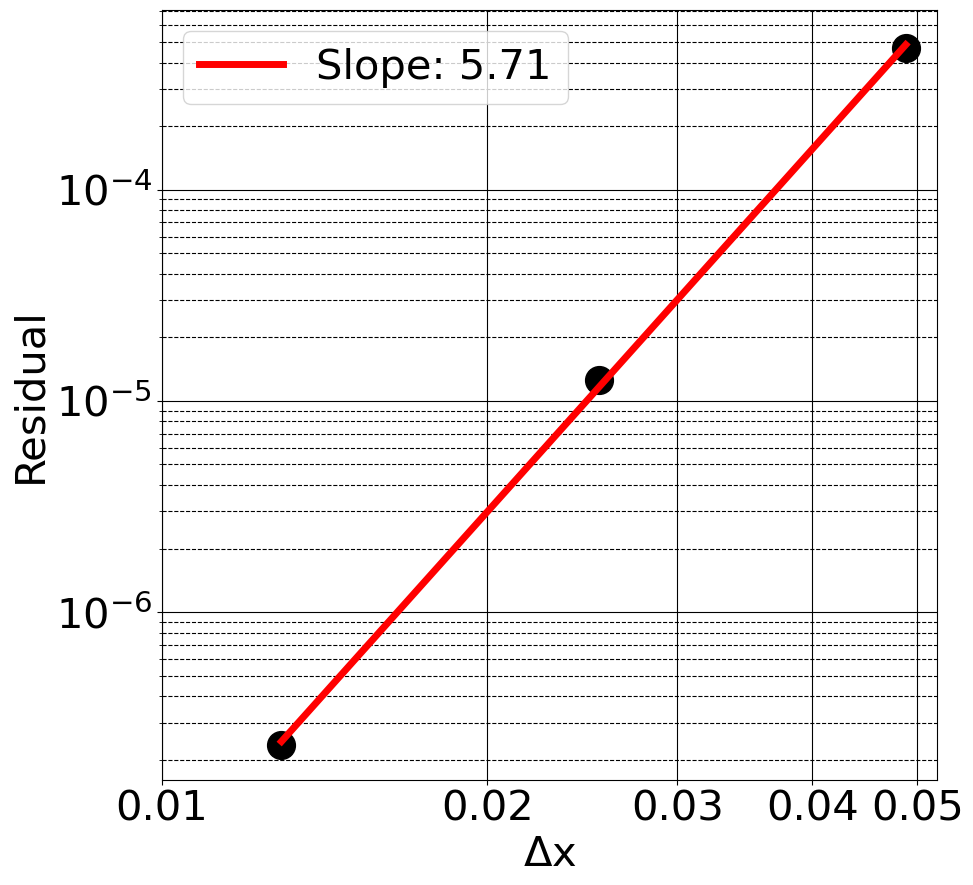}
		\caption{Degree of Appended Polynomial: 6}
	\end{subfigure}
	\caption{Convergence Rate of Final Residuals}
	\label{Fig:Meshless: double circ rect rate of conv residual}
\end{figure}
For this case as well, we found that the residuals of the SOR iterations at the end of each timestep get saturated at values based on the grid resolution and polynomial degree. These saturated values are plotted in \cref{Fig:Meshless: double circ rect rate of conv residual} with average $\Delta x$ for polynomial degrees of 3 and 6. The best fit lines show that the rates of convergence improve with polynomial degree.

\begin{figure}[H]
	\centering
	\includegraphics[width=0.8\textwidth]{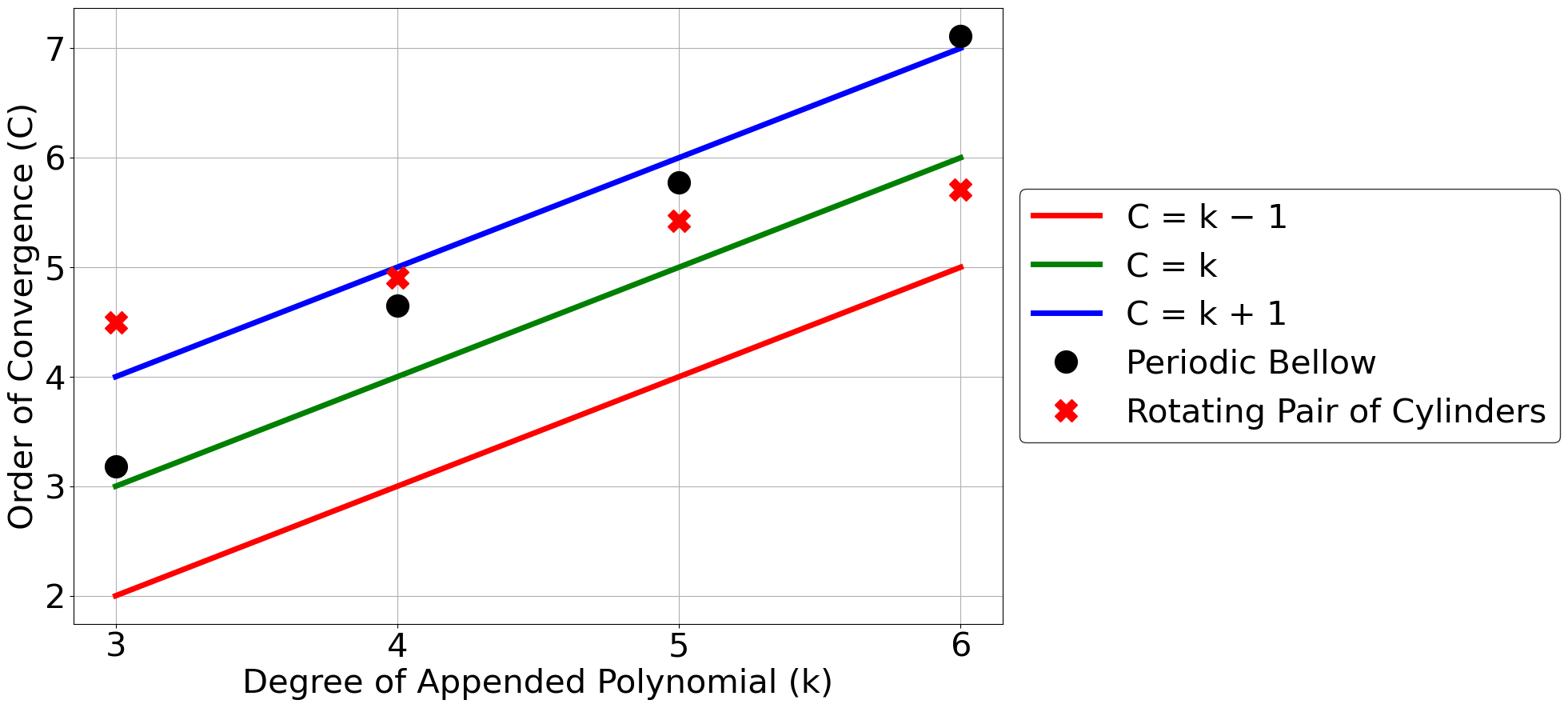}
	\caption{Composite Convergence Rate of Final Residuals with Degrees of Appended Polynomial}
	\label{Fig:Meshless: applications rate of conv residual with polydeg}
\end{figure}
\Cref{Fig:Meshless: applications rate of conv residual with polydeg} is a combined plot of the rates of convergence of the final residuals for the flow in periodic converging-diverging channel and the flow induced by a rotating pair of cylinders. We see that the rate increases monotonically with polynomial degree and follows the reference lines. This is consistent with the observations in the solutions of Poisson equation and fluid flow with analytical solutions.

\section{Summary and Conclusions}
In this paper, we have studied the consistency and accuracy of a high order accurate meshless Navier-Stokes algorithm for incompressible flows in complex domains with all velocity boundary conditions. For incompressible flows with all velocity boundary conditions, the pressure Poisson equation needs to be solved with all Neumann boundary conditions. This results in an ill-conditioned equation for which solutions can be obtained only to an arbitrary level. When solved with finite difference methods, the iterative solver is seen to reach a stationary level of residuals, which is reflective of the inconsistency in the discretization of the derivatives and the source terms. In this paper, we observe that such inconsistency errors, as indicated by the stationary values of the residuals, decrease with grid resolution and order of discretization. Hence, we are of the opinion that it is not necessary to regularize the discrete Poisson equation either by fixing a mean value or the pressure at a reference point. We observe that the convergence is fast (even using an inexpensive SOR solver) without the use of regularization although the iterations reach a stationary level. We further observe that the rate of convergence of these stationary residuals is similar or better than the rate of convergence of the discretization errors. Hence the inconsistencies in the pressure Poisson equation do not appreciably worsen the solution errors in the flow fields.
\par We observe the same in two model flows with analytical solutions, and two complex flows. For all cases, the inter-point spacing and the degree of the appended polynomial are systematically varied. As with the Poisson equation with manufactured solution, we see that the Poisson equation for an actual flow also reaches a stationary level (reflecting the inconsistency), but the inconsistency decreases fast as expected by the high order accuracy of the meshless method. We therefore conclude that for all velocity boundary conditions, one can stop the pressure iterations when they reach the asymptotic level, and not employ regularization conditions such as fixing the mean or a point pressure. We have used here a simple single-grid SOR scheme for its simplicity and economy. However, it can be made to converge much faster by using a multilevel strategy, as in our recent work \cite{radhakrishnan2021non}.

\bibliography{References}

\end{document}